\documentclass[11pt]{amsart}
\pdfoutput=1
\usepackage{amssymb}
\usepackage{latexsym}
\usepackage{graphics}
\usepackage{graphicx}
\usepackage{epstopdf}
\usepackage[all]{xy}
\usepackage{verbatim}
\input epsf

\newtheorem{thm}{Theorem}
\numberwithin{thm}{section}
\newtheorem{prop}[thm]{Proposition}
\newtheorem{definition}[thm]{Definition}
\newtheorem{cor}[thm]{Corollary}
\newtheorem{lemma}[thm]{Lemma}

\newcommand{\startproof}{\textit{Proof:}\ }
\newcommand{\finishproof}{\hfill$\square$}

\def\R{{\mathbb R}}

\begin{document}
\title{Knotting Probability of Equilateral Hexagons}
\author[Hake]{Kathleen Hake}



\address{Department of Mathematics and Statistics, Carleton College, Northfield, MN, 55057}

\date \today

\begin{abstract}

For a positive integer $n\ge 3$, the collection of $n$-sided polygons embedded in $3$-space defines the space of geometric knots. We will consider the subspace of equilateral knots, consisting of embedded $n$-sided polygons with unit length edges. Paths in this space determine isotopies of polygons, so path-components correspond to equilateral knot types. When $n\le 5$, the space of equilateral knots is connected. Therefore, we examine the space of equilateral hexagons. Using techniques from symplectic geometry, we can parametrize the space of equilateral hexagons with a set of measure preserving action-angle coordinates. With this coordinate system, we provide new bounds on the knotting probability of equilateral hexagons.

\end{abstract}

\maketitle

\section{Introduction}

Classically a knot can be defined as a closed, non self-intersecting smooth curve embedded in Euclidean $3$-space. Two knots are considered to be equivalent if one can be smoothly deformed into another. The question of whether or not two given knots are equivalent proves to be a difficult problem. Much of theory is devoted to developing techniques to answer this question. The study of the invariance of knots has been of interest to not only mathematicians but also biologist, physicists, and computer scientists.  Prominent examples of knotting appear in polymers, specifically DNA and proteins. In the early 1970's, it was discovered that enzymes called topoisomerases causes the DNA to change its form. Type II topoisomerases bind to two segments of double-standed DNA, split one of the segments, transport the other through the break, and reseal the break. These studies suggest that the topological configuration, or the knotting, plays a role in understanding the behavior of these enzymes. Sometimes the arbitrary flexibility and lack of thickness in the classical theory of knots does not accurately depict the physical constraints of objects in nature. This inspires questions in the field of physical knot theory and models that seek to capture some of the physical properties.

A question of focus in this paper is about the statistical distribution of knot types as a function of the length. For example, what is the probability that an $n$ edged polygon is knotted? The study of knots from a probabilistic viewpoint provides an understanding of typical knots. There are many ways to model random knots. The model we will consider is that of closed polygonal curves in $\mathbb{R}^3$. A knot is realized by joining $n$ line segments. In addition, we restrict the length of the segments to be equal. We identify each $n$-sided polygonal curve with the $3n$-tuple of vertex coordinates which define it. This gives a correspondence between points in $\mathbb{R}^{3n}$ and $n$-sided polygons in $\mathbb{R}^3$. We consider the $2(n-3)$ dimensional subspace of equilateral polygons equivalent up to translations and rotations. Using techniques from symplectic geometry, we study the space of equilateral hexagons. Suppose $P$ is an equilateral hexagon. We prove that the probability $P$ is knotted is at most $\frac{14-3\pi}{192}<\frac{1}{42}$.


\section{Background}

There are various ways to define a knot, all of which capture the intuitive notion of a knotted loop. We will start by defining a polygonal knot. For any two distinct points in $3$-space, $p$ and $q$, let $[p,q]$ denote the line segment joining them. For an ordered set of distinct points, $(p_1, p_2, \dots, p_n)$, the union of the segments $[p_1, p_2], [p_2,p_3],\dots, [p_{n-1},p_n],$ and $[p_n,p_1]$ is called a closed polygonal curve. If each segment intersects exactly two other segments, intersecting each only at an endpoint, then the curve is said to be simple.

\begin{definition}
	A polygonal knot, $K$, is a simple, closed polygonal curve in $\mathbb{R}^3$.
\end{definition}
 
 If a polygonal knot, $P$, has $n$ vertices, we will call $P$ a $n$-sided polygonal knot. We label the vertices of $P$ as $v_{1}, v_{2}, \ldots, v_{n}$. We call the segments $[v_i,v_{i+1}]$ the edges of $P$ and label the edges of $P$ as $e_{1}, e_{2}, \ldots, e_{n}$, where $e_{1}=[v_1,v_2], e_2=[v_2, v_3], \dots, e_{n-1}=[v_{n-1},v_n]$, and $e_n=[v_n,v_1]$. In addition, we will select a distinguished vertex, $v_{1}$, called a root and a choice of orientation.

\begin{figure}[h]
\centerline{\includegraphics[width=.35\textwidth]{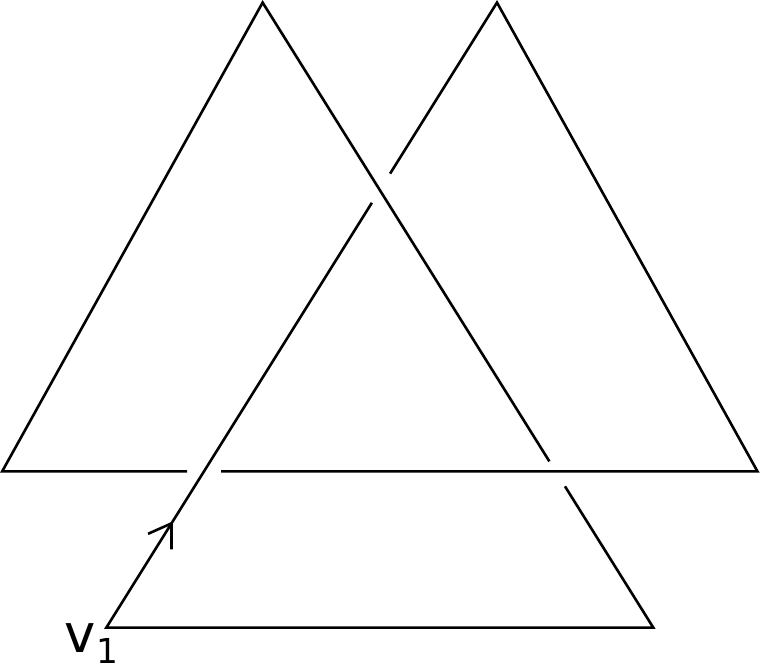}\hspace{1.5cm}\includegraphics[width=.3\textwidth]{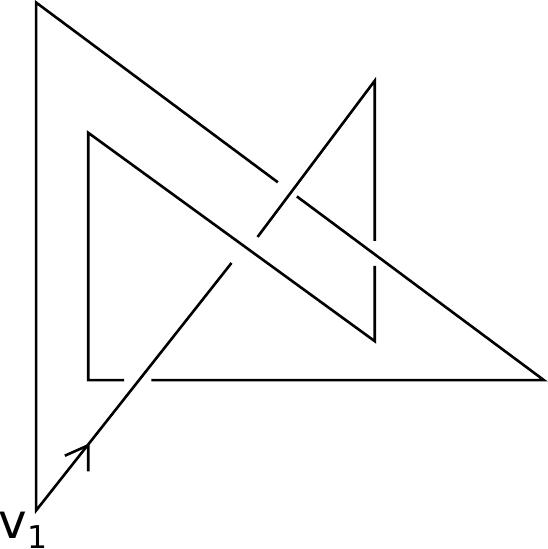}}
\caption{The figure on the left shows a 6-sided, rooted, oriented polygonal trefoil knot. The figure on the right shows a 7-sided, rooted, oriented polygonal figure-8 knot.}
\label{fig:label}
\end{figure}

 With a distinguished vertex and orientation, we can view $P$ as a point in $\mathbb{R}^{3n}$ by listing the coordinates of the vertices starting with $v_1$ then following the orientation. Not all points in $\mathbb{R}^{3n}$ will correspond to simple polygonal curves. Therefore we define the discriminant set, in the spirit of Vassiliev \cite{Birman1993}.

\begin{definition} The discriminant, $\Sigma^{(n)}$, is all points in $\mathbb{R}^{3n}$ that correspond to non-embedded polygonal knots. 
\end{definition}
A polygonal knot in $\R^{3}$ fails to be embedded when two or more of the edges intersect. For an $n$-sided polygonal knot there are $\frac{1}{2}n(n-3)$ pairs of non-adjacent edges. So $\Sigma^{(n)}$ is the union of the closure of the $\frac{1}{2}n(n-3)$ real semi-algebraic cubic varieties, each piece consisting of polygons with a single intersection between non-adjacent edges\cite{Randell1}\cite{Randell2}. By excluding these singular points, we are left with an open set in $\mathbb{R}^{3n}$ corresponding to embedded polygons in $\mathbb{R}^3$.
\begin{definition}
The embedding space for rooted, oriented $n$-sided polygonal knots, denoted $Geo(n)$, is defined to be $ \mathbb{R}^{3n}-\Sigma^{(n)}$.
\end{definition}
 The space $Geo(n)$ is an open $3n$-dimensional manifold. A continuous path $h:[0,1]\to Geo(n)$ is an isotopy of polygonal knots. 
 
\begin{definition} Two $n$-sided polygonal knots are geometrically equivalent if they lie in the same path-component of $Geo(n)$.
\end{definition}

Path components are in bijective correspondence with the geometric knot types realizable with $n$ edges. Any polygon that is in the same path-component of the regular, planar $n$-gon is then geometrically equivalent to the unknot. 

Next we will consider polygonal knots with unit length edges. Consider the function $f:Geo(n)\to \mathbb{R}^{n}$ where $(v_{1},v_{2},\ldots, v_{n})\mapsto (||v_{1}-v_{2}||,||v_{2}-v_{3}||,\ldots,||v_{n}-v_{1}||)$. 

\begin{definition}
	Let $f^{-1}((1,1,\ldots,1))=Equ(n)$ be the embedding space for rooted, oriented, $n$-sided equilateral knots.
\end{definition}

Since $Equ(n)$ is the preimage of the smooth map $f$ at the regular value $(1,1,\ldots,1)$, $Equ(n)$ is a $2n$-dimensional manifold. Similar to the space of geometric knots, path-components correspond to the equilateral knot types realizable with $n$ edges. In this paper, we will focus on equilateral polygonal knots.

Every triangle is planar. A quadrilateral can be folded along a diagonal to become planar. It is also known that any pentagon can be deformed to a planar pentagon \cite{Randell2}. 
Therefore $Equ(n)$ is connected for $n\le 5$ and the case of hexagons is the first interesting example. Jorge Calvo \cite{Calvo} proves that $Equ(6)$ has five path-components. One component of $Equ(6)$ corresponds to the unknot, two to the right-handed trefoil and two to the left-handed trefoil. In order to distinguish between the different components, he introduces new knot invariants for equilateral hexagonal knots. First let $H=(v_{1},v_{2},\ldots,v_{6})\in Equ(6)$. 

\begin{definition}
Let $H\in Equ(6)$. The curl of $H$, denoted $curl(H)$, is defined by $curl(H)=\text{sign}((v_{3}-v_{1})\times(v_{5}-v_{1})\cdot(v_{2}-v_{1}))$.	
\end{definition}

 If $v_{1}, v_{3}$ and $v_{5}$ are on the $xy$-plane oriented in a counter-clockwise orientation, then $curl(H)$ denotes the sign of the $z$-coordinate of $v_{2}$. So, in a sense, some knots curl up while others curl down. We will describe a second invariant of the hexagonal knot that distinguishes its topological knot type.
 
 \begin{definition}
 	Define $T_i$ to be the interior of the triangular disk spanned by $(v_{i-1}, v_{i}, v_{i+1})$.
 \end{definition} 

\begin{figure}[h]
\centerline{\includegraphics[width=.35\textwidth]{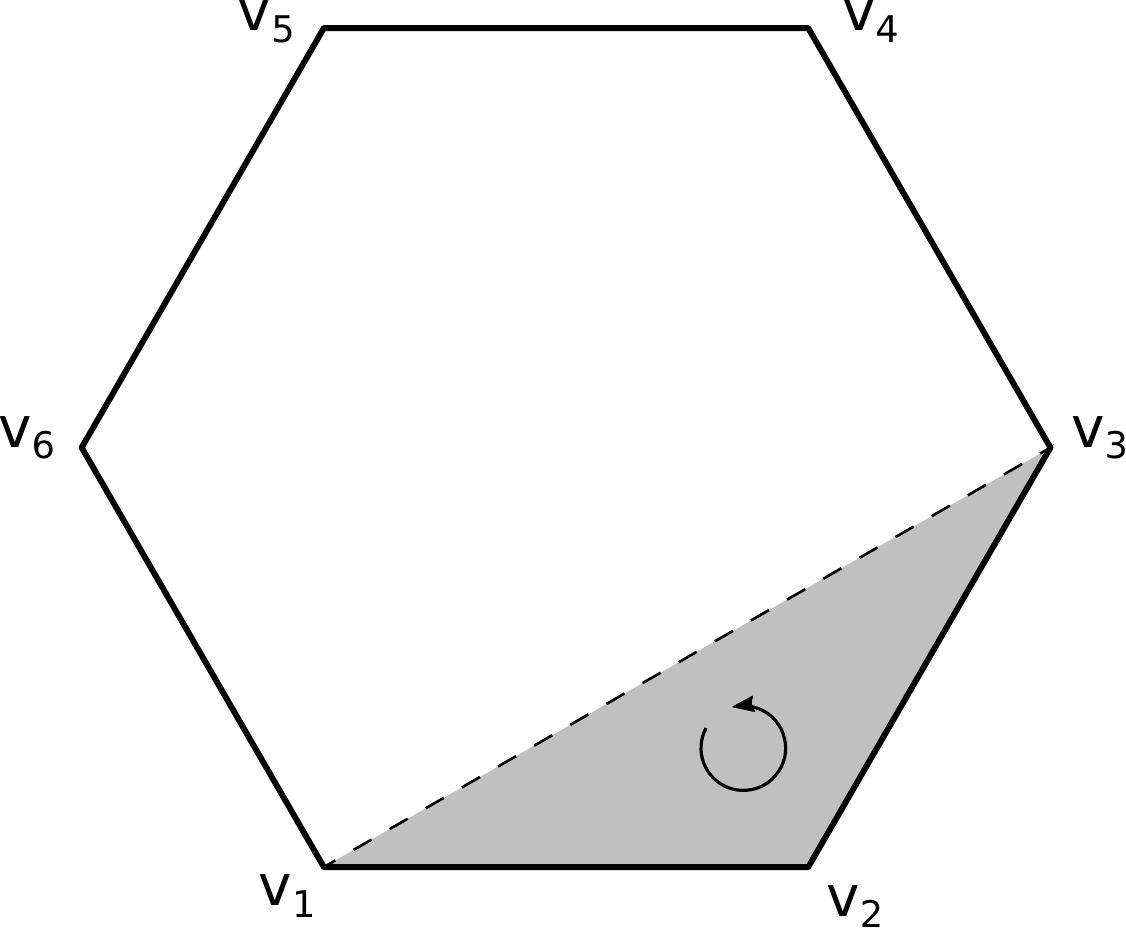}}
\caption{In this figure triangular disk $T_2$ is shaded and the orientation from the right-hand rule is shown.}
\label{fig:T2}
\end{figure}

Using a right-hand rule, we orient each $T_i$ as shown in Figure \ref{fig:T2}.

\begin{definition}
	Define $\Delta_{i}$, for $i=2,4,$ and $6$, to be the algebraic intersection number of T$_{i}$ with $H$.	
\end{definition} 
 
  
  \begin{figure}[h]
\centerline{\includegraphics[width=.35\textwidth]{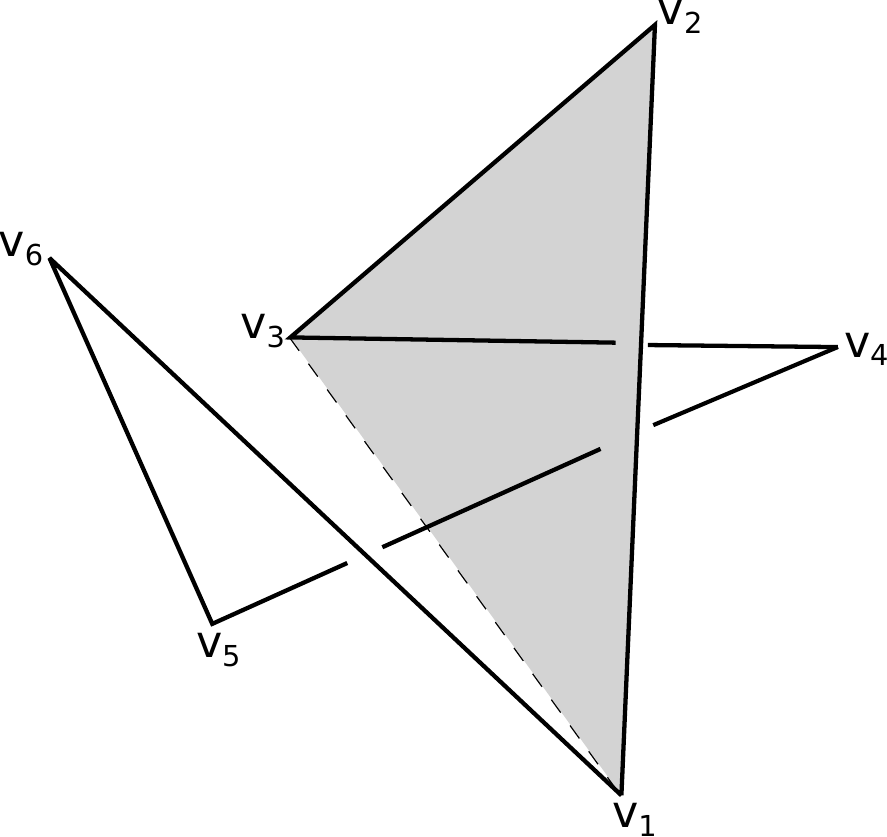}}
\caption{This figure shows an example of a hexagonal knot in which $\Delta_2=1$.}
\label{fig:D2}
\end{figure}

 The following Lemma distinguishes topological knot type using the algebraic intersection numbers.

\begin{lemma}(\cite{Calvo})\label{algintersection}
Let $H\in Equ(6)$. Then
\begin{enumerate}
\item $H$ is a right-handed trefoil iff $\Delta_{i}=1$ for all $i$,
\item $H$ is a left-handed trefoil iff $\Delta_{i}=-1$ for all $i$,
\item $H$ is an unknot iff $\Delta_{i}=0$ for some $i\in\{2,4,6\}$.
\end{enumerate}
\end{lemma}

 Combining the notion of curl with the appropriate intersections from Lemma \ref{algintersection} we arrive at Calvo's Geometric Knot Invariant, Joint Chirality-Curl.
 
 \begin{definition}(\cite{Calvo})
 	Let $H\in Equ(6)$. Define Joint Chirality-Curl
 	$$J(H)=(\Delta_{2}\Delta_{4}\Delta_{6},\Delta_{2}^2\Delta_{4}^2\Delta_{6}^2 curl(H)).$$
 	
 \end{definition} 
 
 The Joint Chirality-Curl distinguishes between the five components of $Equ(6)$.
 
 \begin{thm}(\cite{Calvo}): Let $H\in Equ(6)$. Then

\[
J(H) =
\begin{cases}
(0,0) & \text{iff } $H$\text{ is unknot}\\
(1,c) & \text{iff } $H$ \text{ is right-trefoil with } $curl(H)=c$\\
(-1,c)& \text{iff } $H$\text{ is left-trefoil with }$curl(H)=c$
\end{cases}
\]
\end{thm}

The five components of $Equ(6)$ are due to the choice of a root and orientation. Consider the automorphisms $r$ and $s$ on $Equ(6)$ defined by 
$$ r\langle v_1, v_2, v_3, v_4, v_5, v_6\rangle=\langle v_1, v_6, v_5, v_4, v_3, v_2\rangle$$
$$ s\langle v_1, v_2, v_3, v_4, v_5, v_6\rangle=\langle v_2, v_3, v_4, v_5, v_6, v_1\rangle.$$

These automorphisms act on $Equ(6)$ by reversing or shifting the order of the vertices of each hexagon. They generate the dihedral group of order twelve. 

\begin{thm}(\cite{Calvo})\label{2.5Calvo}
 Suppose $\Gamma$ is a subgroup of the dihedral group $\langle r,s\rangle$. Then $Geo(6)/\Gamma$ has five components if and only if $\Gamma$ is contained in the index-$2$ subgroup $\langle s^2,rs \rangle$. Otherwise, $Geo(6)/\Gamma$ has three components.
\end{thm}

\begin{cor}(\cite{Calvo}) The spaces $Geo(6)/\langle s\rangle $  of non-rooted oriented embedded hexagons, and $Geo(6)/\langle r,s\rangle$ of non-rooted non-oriented embedded hexagons, each consist of three path-components.
\end{cor}

Next we will discuss definitions and results from symplectic geometry, specifically toric symplectic manifolds \cite{symplectic}, that apply to knot spaces.

\begin{definition}
A symplectic manifold, $M$, is an even dimensional manifold with a closed, non-degenerate $2$-form, $\omega$, called the symplectic form. 
\end{definition}

Since $\omega$ is non-degenerate, there is a canonical isomorphism between the tangent and cotangent bundles, namely 

$$TM\mapsto T^*M: X\to \iota(X)\omega=\omega(X,\cdot ).$$ 

\begin{definition}
A symplectomorphism of a symplectic manifold $(M,\omega)$ is a diffeomorphism $\psi\in Diff(M)$ that preserves the symplectic form. The group of symplectomorphisms of $M$ is denoted $Symp(M,\omega)$.
\end{definition}

Since $\omega$ is nondegenerate the homomorphism $T_q M\to T_q^*M: v\mapsto \iota(v)\omega_q$ is bijective. Thus there is a one-to-one correspondence between vector fields and $1$-forms via
$$\chi (M)\to \Omega^1(M):X\mapsto \iota (X)\omega.$$ 

\begin{definition}
	A vector field $X\in \chi (M)$ is called symplectic if $\iota (X)\omega$ is closed. Denote the space of symplectic vector fields by $\chi (M,\omega)$.
\end{definition} 

\begin{prop}\cite{symplectic}\label{symp}
	Let $M$ be a closed manifold. If $t\mapsto \psi_t\in Diff(M)$ is a smooth family of diffeomorphims generated by a family of vector fields $X_t\in \chi(M)$ via
	$$ \frac{d}{dt}\psi_t=X_t \circ \psi_t, \hspace{1cm} \psi_0=id,$$
	then $\psi_t\in Symp(M,\omega)$ for every $t$ if and only if $X_t\in \chi(M,\omega)$ for every $t$.
\end{prop}

Now consider a smooth function $H:M\to \mathbb{R}$. 

\begin{definition}
	The vector field $X_H:M\to TM$ determined by identity $dH=\iota(X_H)\omega$ is called the Hamiltonian vector field associated to $H$.
\end{definition}

If $M$ is closed, then by Proposition \ref{symp}, the vector field $X_H$ generates a smooth $1$-parameter group of diffeomorphisms $\phi^t_H\in Diff(M)$ such that $$\frac{d}{dt} \phi^t_H=X_H\circ \phi^t_H,\hspace{1cm} \phi^0_H=id,$$
called a Hamiltonian flow associated to $H$. The identity $$dH(X_H)=(\iota(X_H)\omega)(X_H)=\omega(X_H,X_H)=0$$ shows that $X_H$ is tangent to level sets.

 A useful example to consider is the unit sphere $S^2$, where $\omega$ is the standard area form. If $S^2=\{(x_1, x_2, x_3) : \sum_j x_j^2=1 \}$,  then $\omega_x(u,v)= \langle x,u\times v\rangle $ for $u,v\in T_x S^2$. Consider cylindrical polar coordinates $(\theta,x_3)$ for $\theta\in [0,2\pi)$ and $x_3 \in [-1,1]$. Let $H$ be the height function $x_3$ on $S^2$. The level sets are circles at constant height. The Hamiltonian flow $\phi^t_H$ rotates each circle at constant speed and $X_H$ is the vector field $\frac{\partial }{\partial \theta}$. Thus $\phi^t(H)$ is the rotation of the sphere about its vertical axis through the angle $t$.

Consider a smooth map $[0,1]\times M\to M:(t,q)\mapsto \psi_t(q)$ such that $\psi_t\in Symp(M,\omega)$ and $\psi_0=id$. A family of such symplectomorphisms is called a symplectic isotopy of $M$. The isotopy is generated by a unique family of vector fields $X_t:M\to TM$ such that $\frac{d}{dt}\psi_t=X_t\circ\psi_t$. If all of the $1$-forms are exact then there exists a smooth family of Hamiltonian functions $H_t:M\to \mathbb{R}$ such that $\iota (X_t)\omega=dH_t$. In this case, $\psi_t$ is called a Hamiltonian isotopy.

\begin{definition}
	A symplectomorphism, $\psi$, is called Hamiltonian if there exists a Hamiltonian isotopy $\psi_t\in Symp(M,\omega)$ from $\psi_0=id$ to $\psi_1=\psi$.
\end{definition}

\begin{definition}
	A Hamiltonian action of $S^1$ on $(M,\omega)$ is a $1$-parameter subgroup $\mathbb{R}\to Symp(M): t\mapsto \psi_t$ of $Symp(M)$ where $\psi_t=id$ and which is the integral of a Hamiltonian vector field $X_H$.
\end{definition}

The Hamiltonian function $H:M\to \mathbb{R}$ in this case is called the moment map. If $k$ such symmetries commute we have an action of a torus, $T^{k}$, on $M$. Then the moment map, $\mu: M\to\R^{k}$ yields a $k$-dimensional vector of conserved quantities. If $k$ is half the dimension of $M$, then $M$ is called toric symplectic. From theorems of Atiyah\cite{Atiyah:1982re} and Guillemin-Sternberg\cite{GuillStern}, the image of $\mu$ is a convex polytope, $P$, called the moment polytope. Moreover, the vertices of the moment polytope are the images under $\mu$ of the fixed point of the Hamiltonian torus action. In addition, the torus action preserves the fibers of the moment map. If we can invert $\mu$, we get a map $\alpha:P\times T^{n}\to M$ called the action-angle map. 

The previous example of the unit sphere $S^2$ is a toric symplectic manifold, with circle action rotation about the $z$-axis. The moment map $H:S^2\to S^1$ is the height function, the conserved quantity as the sphere rotates. The image of $H$ is a convex polytope, namely the interval $[-1,1]$. The fibers of $H$ are horizontal circles of constant height, which are preserved under the action. Lastly the circle, $S^1$, is half the dimension of $S^2$.

The toric symplectic structure on the sphere naturally carries over to a toric symplectic structure on the product of spheres. This gives a toric symplectic structure on the space of open random walks or open polygons. We will consider the subspace of closed random walks. Let $Pol(n)$ be the $2n$ dimensional space of possibly singular polygons in $\mathbb{R}^3$ with edgelengths one. We will consider the quotient space $Pol_0(n)=Pol(n)/\text{SO}(3)$ of equilateral polygons up to translations and rotations. Jason Cantarella and Clayton Shonkwiler \cite{JC} describe the almost toric symplectic structure of $Pol_0(n)$. We summarize some of the important information below. To define the toric action, consider any triangulation, $T$, of an equilateral planar regular $n$-gon. Let $d_{i}$ be the lengths of the $n-3$ diagonals of the triangulation. These diagonals, along with the edges on the polygon, form $n-2$ triangles which each obey $3$ triangle inequalities. Therefore the lengths of the diagonals and the edge lengths must obey a set of $3(n-2)$ triangle inequalities, called the triangulation inequalities. 

\begin{thm}\cite{Kapovich}\cite{Howardmanon}\cite{Hitchin}
	The following facts are known:
	\begin{itemize}
		\item $Pol_0(n)$ is a possibly singular $(2n-6)$-dimensional symplectic manifold. The symplectic volume is equal to the standard measure.
		\item To any triangulation $T$ of the standard $n$-gon we can associate a Hamiltonian action of the torus $T^{n-3}$ on $Pol_0(n)$, where $\theta_i$ acts by folding the polygon around the $i^{th}$ diagonal of the triangulation.
		\item The moment map $\mu: Pol_0(n)\mapsto \mathbb{R}^{n-3}$ for a triangulation $T$ records the lengths $d_i$ of the $n-3$ diagonals of the triangulation.
		\item The inverse image $\mu^{-1}(int(P))\subset Pol_0(n)$ of the interior of the moment polytope $P$ is an toric symplectic manifold.
	\end{itemize}
\end{thm}

\begin{figure}[h]
\centerline{\includegraphics[width=.5\textwidth]{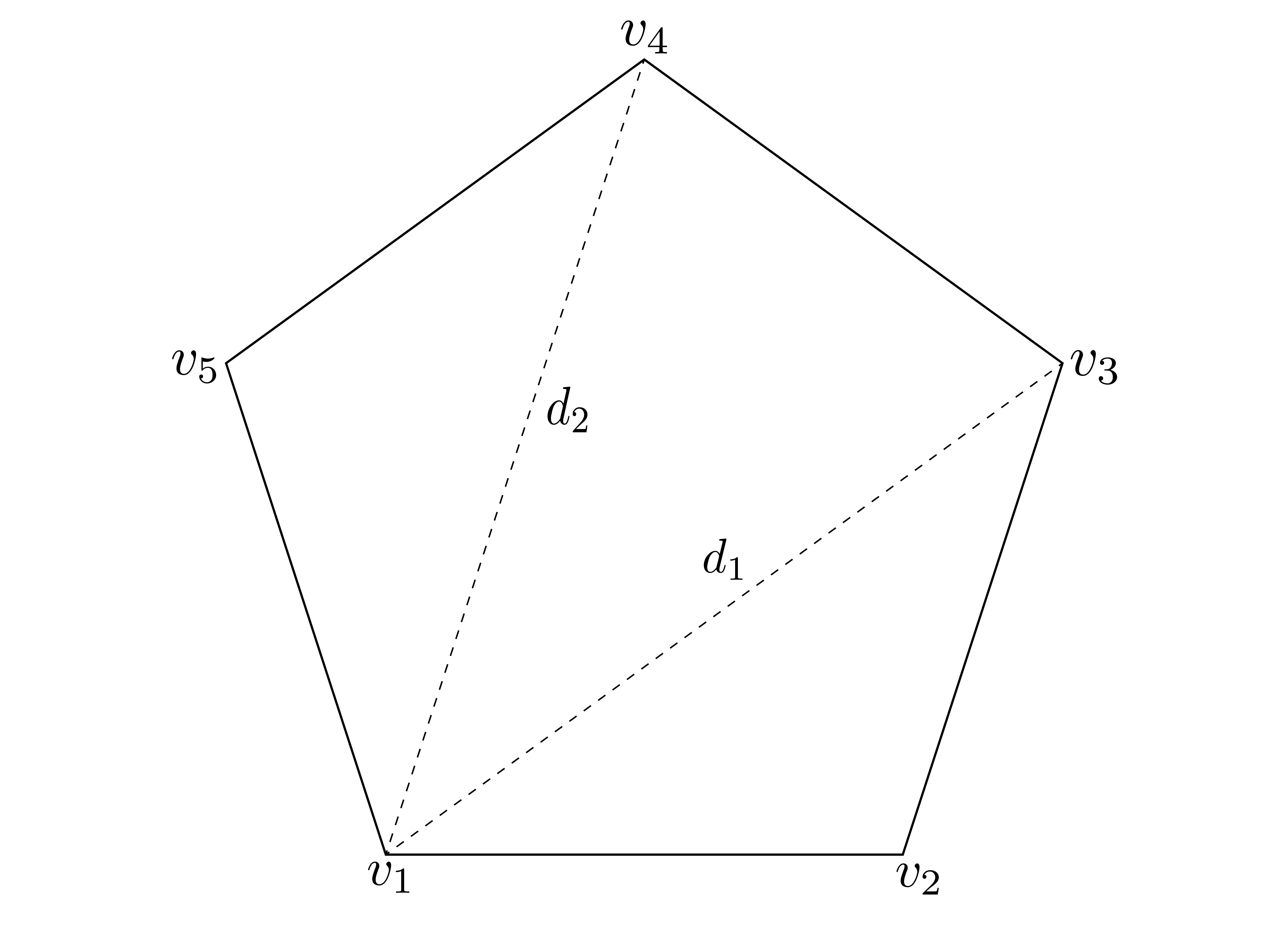}\includegraphics[width=.55\textwidth]{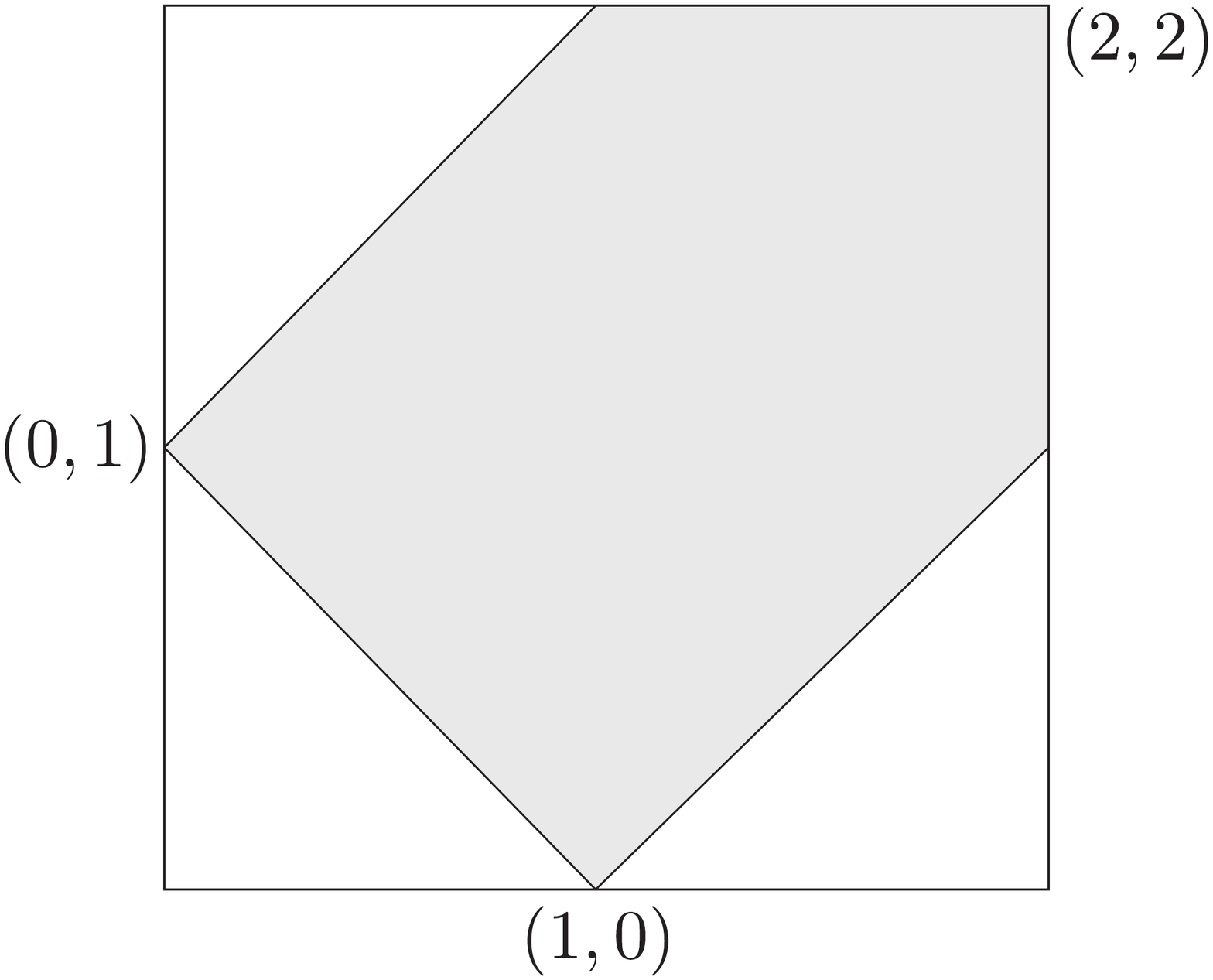}}
\caption{The left image shows the fan triangulation of an equilateral pentagon, where all diagonals share a common vertex. The lengths of the diagonals, $d_1$ and $d_2$, satisfy six triangle inequalities. The figure on the left shows the moment polytope of $Pol_0(5)$ corresponding to the fan triangulation.}
\label{pentagon}
\end{figure}

The moment polytope, $P_n$, is defined by the triangulation inequalities for $T$. The vertices of the moment polytope represent degenerate polygons which extremize several triangulation inequalities. Figure \ref{pentagon} shows a triangulation of an equilateral pentagon and the corresponding moment polytope.

\begin{figure}[h]
\centerline{\includegraphics[width=.35\textwidth]{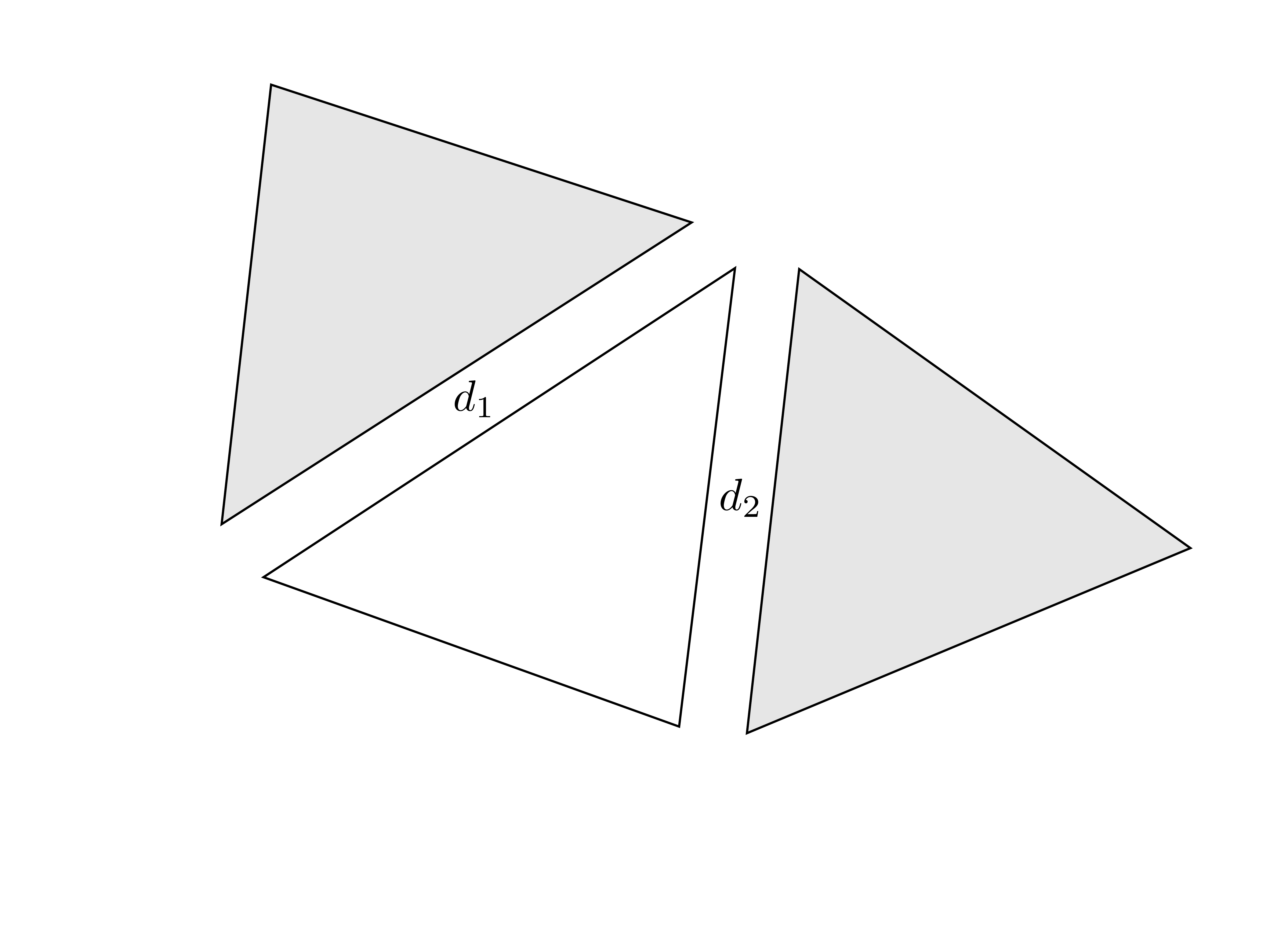}\includegraphics[width=.35\textwidth]{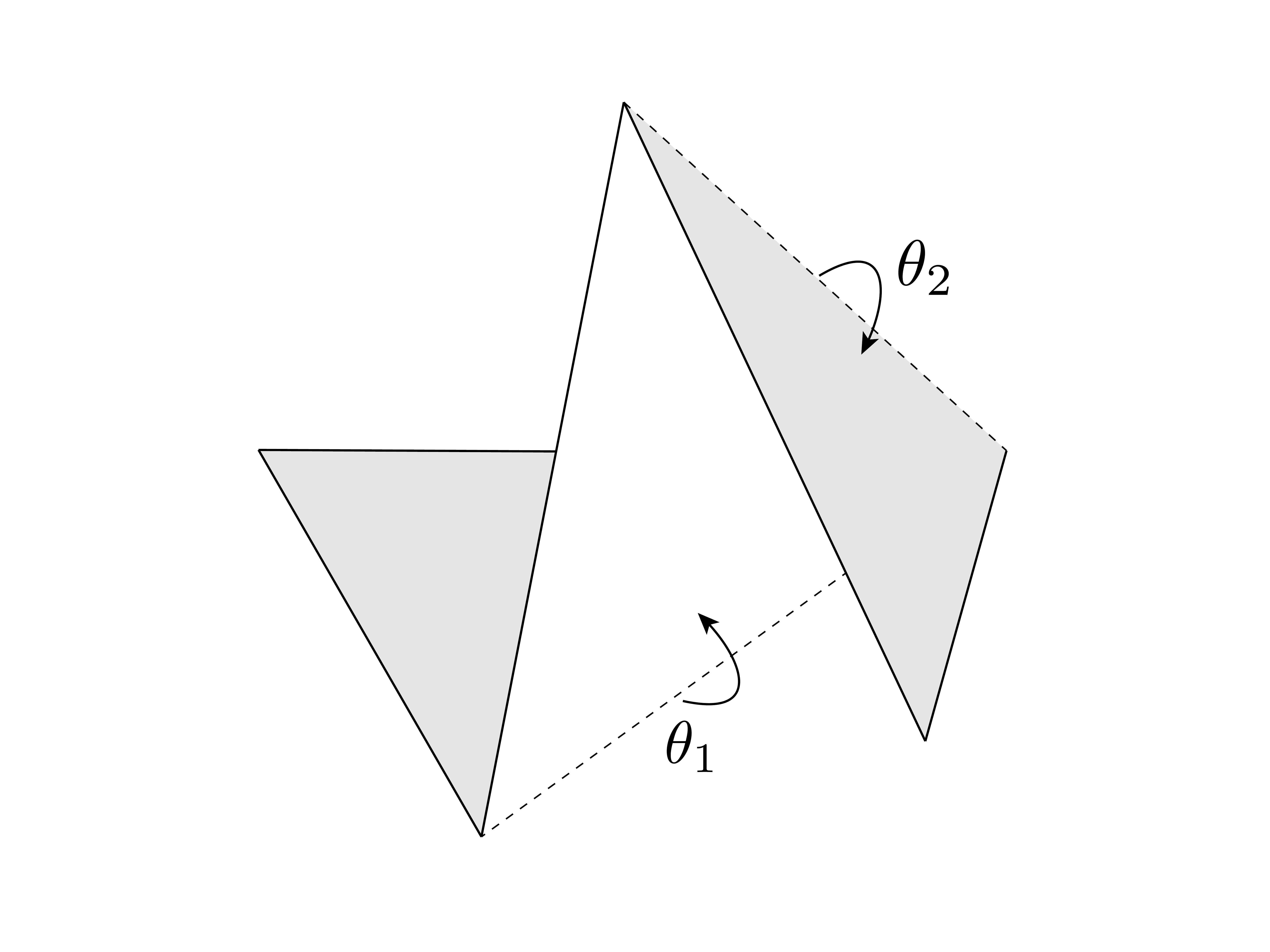}\includegraphics[width=.35\textwidth]{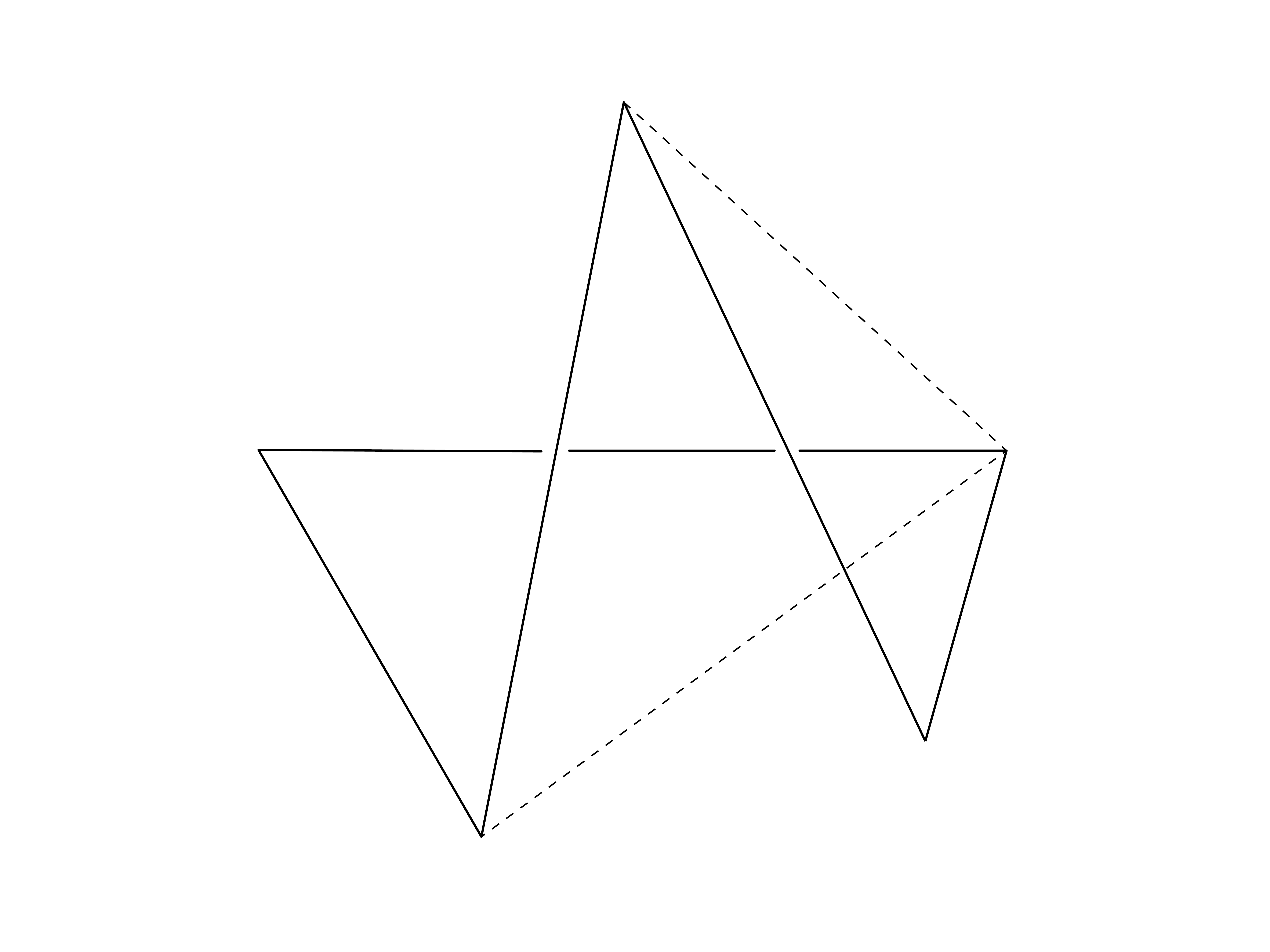}}
\caption{The figure shows how to construct an equilateral pentagon from the action-angle map $\alpha: P_5\times T^2\mapsto Pol_0(5)$ for the fan triangulation. A point $(d_1, d_2)$ in the moment polytope gives the information needed to construct three triangles. Then a point $(\theta_1, \theta_2)\in T^2$ gives instruction on how to attach the triangles along the diagonals. The boundary of the triangulated surface is the equilateral pentagon.}
\label{pentagon2}
\end{figure}

The action-angle map $\alpha:P_n\times T^{n-3} \mapsto Pol_0(n)$ for a triangulation $T$ is given by first constructing the $n-2$ triangles using the diagonal lengths, $d_i$, and edge lengths of $1$ and then joining them in $3$-space with dihedral angles given by the $\theta_i$. The polygon is the boundary of this triangulated surface. This construction only makes sense for polygons equivalent up to translations and rotations, which is why the quotient by $\text{SO}(3)$ is necessary. An example of an equilateral pentagon is shown in Figure \ref{pentagon2}.

The following theorem will be used in Section 4 when calculating the knotting probability of equilateral hexagons. 

\begin{thm}(Duistermaat-Heckman)\cite{DuistHeck}\label{DH}
Suppose $M$ is a $2n$-dimensional toric symplectic manifold with moment polytope $P$, $T^{n}$ is the $n$-torus and $\alpha$ inverts the moment map. If we take the standard measure on the $n$-torus and the uniform measure on $\text{int}(P)$, then the map $\alpha:\text{int}(P)\times T^{n}\to M$ parametrizing a full-measure subset of $M$ in action-angles coordinates is measure-preserving. In particular, if $f:M\to \mathbb{R}$ is any integrable function then $$\int_M f(x)\text{ d}m=\int_{P\times T^n} f(d_1,\cdots,d_n,\theta_1,\cdots,\theta_n)\text{ dVol}_{\mathbb{R}^n}\wedge d\theta_1\wedge\cdots \wedge d\theta_n $$ and if $f(d_1, \cdots,d_n,\theta_1,\cdots,\theta_n)=f_d(d_1, \cdots,d_n)f_\theta(\theta_1,\cdots,\theta_n)$ then $$\int_M f(x)\text{ d}m=\int_{P} f_d(d_1,\cdots,d_n)dVol_{\mathbb{R}^n}\int_{T^n} f_\theta(\theta_1,\cdots,\theta_n)d\theta_1\wedge\cdots \wedge d\theta_n. $$\\

\end{thm}


\section{Symplectic Structure of the space of equilateral hexagons}

\subsection{Action-Angle Coordinates}
 In order to describe the action-angle coordinates on the space of equilateral hexagons, we first must consider the quotient space of $Equ(6)$.
 \begin{definition}
 	Let $Equ_0(6)=Equ(6)/\text{SO}(3)\times\mathbb{R}^3$ be the embedding space of rooted, oriented equilateral hexagonal knots up to translations and rotations.
 \end{definition}

 Let $H=(v_{1},v_{2},v_{3},v_{4},v_{5},v_{6})\in Equ(6)$. We can translate $H$ so that   $v_{1}=(0,0,0)$. Additionally we rotate $H$ so that $v_{3}$ is on the positive $x$-axis and $v_5$ on the upper-half $xy$-plane.  Therefore $v_1, v_3$, and $v_5$ are on the $xy$-plane in a counter-clockwise orientation. In this section, we will consider this to be the standard position for $H\in Equ_0(6)$.
 
 Next we can choose any triangulation of the standard planar equilateral hexagon to form our action-angle coordinates. We will use one of the triangulations that has a central triangle.
 
 \begin{definition} Let the $T_{135}$ triangulation be the triangulation of the regular planar equilateral hexagon which has diagonals connecting $v_{1}$ to $v_{3}$, $v_{3}$ to $v_{5}$, and $v_{5}$ to $v_{1}$, with lengths $d_{1}$, $d_{2}$, and $d_{3}$, respectively.
 \end{definition}
 
 \begin{figure}[h]
\centerline{\includegraphics[width=.4\textwidth]{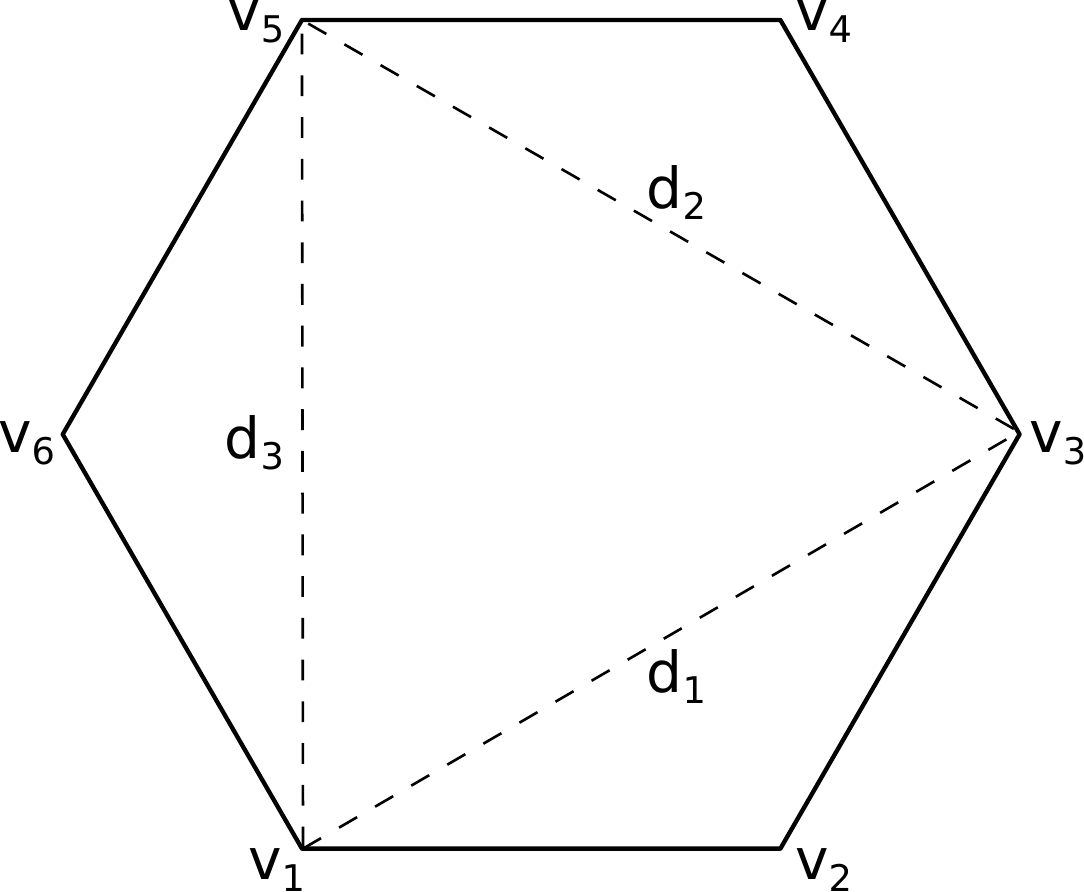}}
\caption{This figure shows the $T_{135}$ triangulation of an equilateral hexagon.}
\label{fig:label}
\end{figure}
 
 The lengths of the diagonals of the $T_{135}$ triangulation obey the following triangulation inequalities:
 \begin{equation*}
\begin{aligned}[c]
0\le& d_1 \le2,\\
0\le& d_2 \le2,\\
0\le& d_3 \le2,\\
\end{aligned}
\quad \text{and}\quad
\begin{aligned}[c]
d_3\le& d_1+d_2,\\
d_1\le&d_3+d_2,\\
d_2\le&d_3+d_1.\\
\end{aligned}
\end{equation*}

\begin{definition}
	The $T_{135}$ triangulation polytope, $P_6$, is the moment polytope for $Pol_0(6)$ corresponding to the $T_{135}$ triangulation and is determined by the triangulation inequalities.
\end{definition}

 \begin{figure}[h]
\centerline{\includegraphics[width=.5\textwidth]{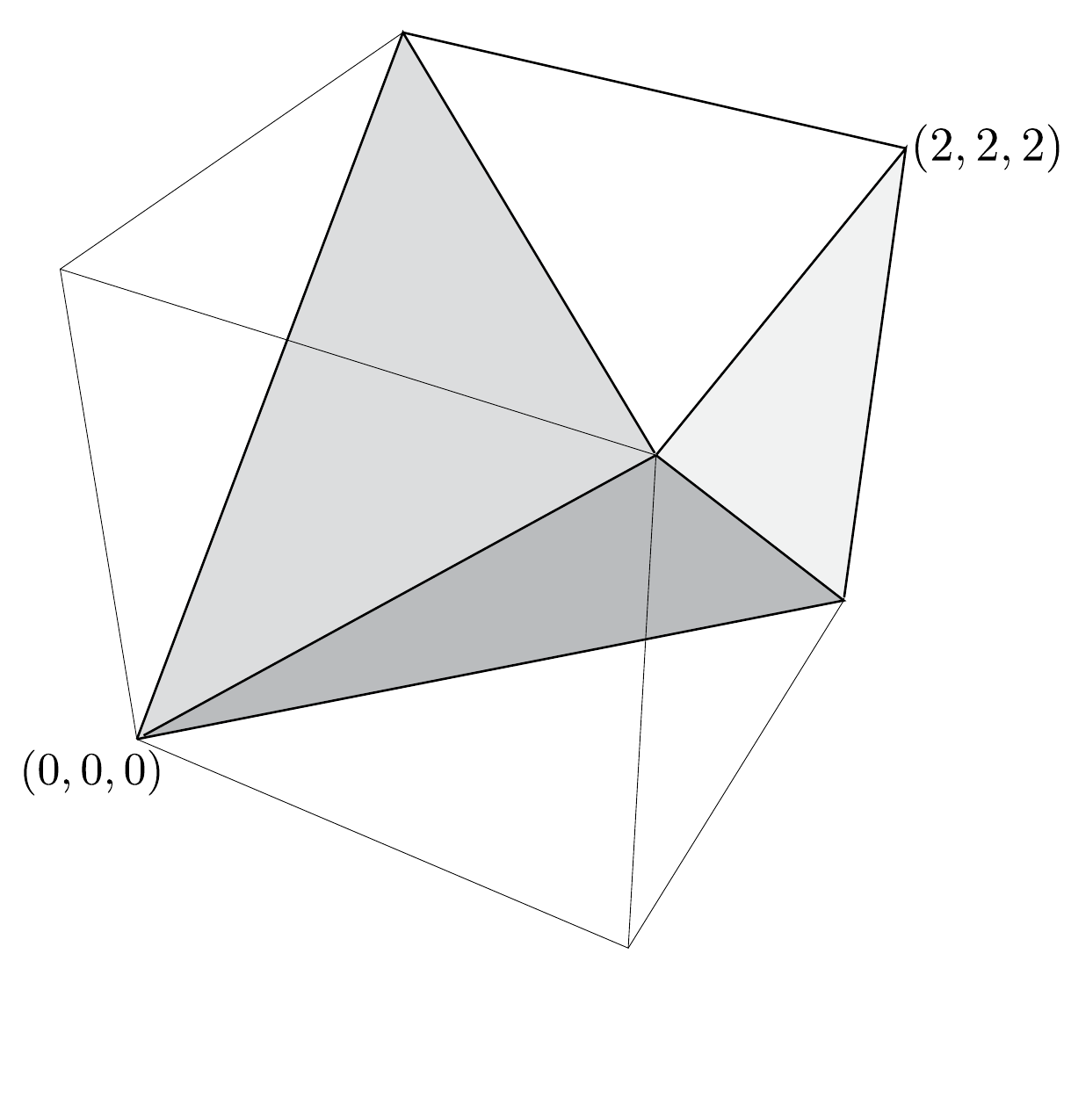}}
\caption{This figure shows the $T_{135}$ triangulation polytope.}
\label{fig:polytope}
\end{figure}

Let $\theta_{i}$ be the dihedral angle around diagonal $d_{i}$, where the regular planar hexagon has all angles $\pi$. Then the action-angle map for $T_{135}$, $\alpha:P_6\times T^3\mapsto Pol_0(6)$ allows us to parametrize any $H\in Equ_{0}(6)$ as $H=(d_{1},d_{2}, d_{3}, \theta_{1}, \theta_{2}, \theta_{3})$. To construct an equilateral hexagonal knot in $Equ_0(6)$, first choose a point $(d_1, d_2, d_3)\in P_6$ and construct four triangles: one with lengths $d_1, d_2$, and $d_3$ and three isosceles triangles with two side lengths $1$ and third side $d_i$. The triangle with side lengths $d_1, d_2$, and $d_3$ is placed on the $xy$-plane with $v_1$ the origin and $v_3$ on the positive $x$-axis. Then a point $(\theta_1,\theta_2,\theta_3)$ in the torus $T^3$ gives instructions on how to connect the three remaining triangles.

\begin{figure}[h]
\centerline{\includegraphics[width=.6\textwidth]{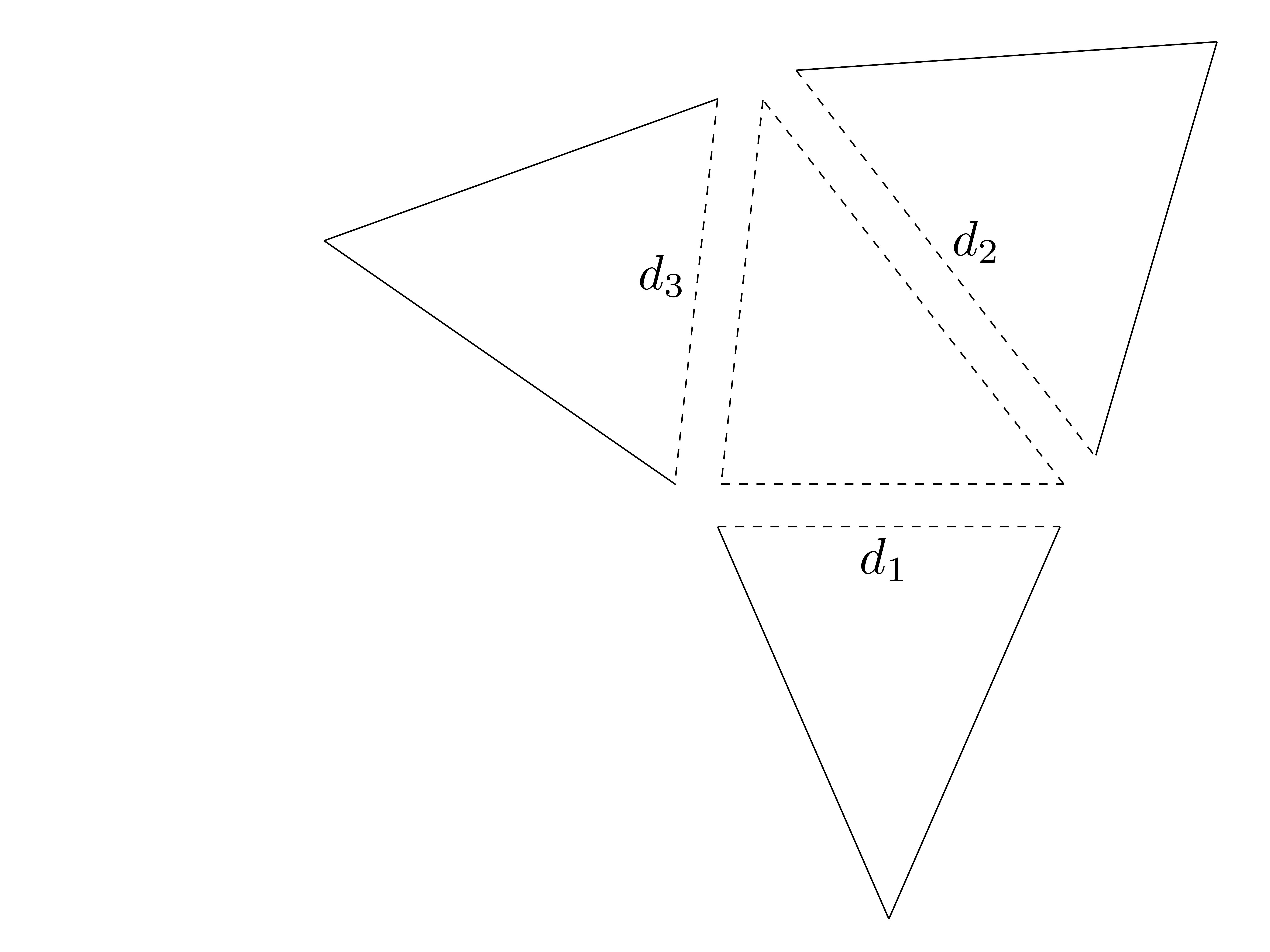}\includegraphics[width=.45\textwidth]{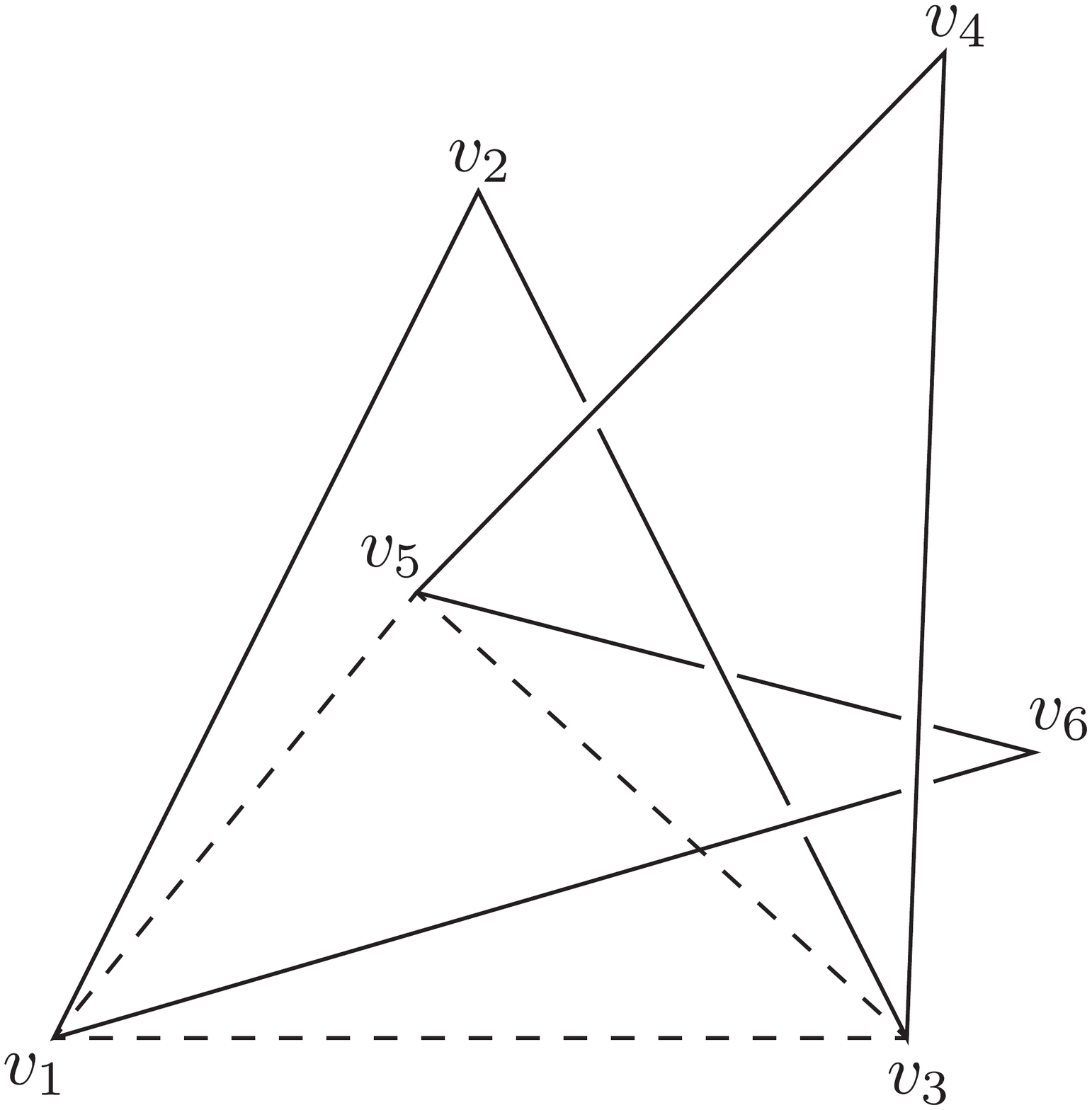}}
\caption{Given a point $(d_1,d_2,d_3)\in P_6$ four triangles are formed. Then given a triple of angles, the triangles are connected to form an equilateral hexagonal trefoil.}
\label{fig:label}
\end{figure}

For $H\in Equ_0(6)$ in standard position, the action-angle coordinates arising from the $T_{135}$ triangulation gives the following coordinates for the vertices of $H$:

\begin{align*}
v_1 =& \Big(0,0,0\Big)\\
v_{2}=& \Big(\frac{d_{1}}{2}, \frac{1}{2}\sqrt{4-(d_{1})^2}\text{ cos}(\theta_{1}),\frac{1}{2}\sqrt{4-(d_{1})^2}\text{ sin}(\theta_{1})\Big),\\ 
v_{3}=& \Big(d_{1},0,0\Big),\\ 
v_{4}=&\Big(\frac{3(d_{1})^2-(d_{2})^2+(d_{3})^2}{4d_{1}}-\frac{d}{4d_{1}d_{2}}\sqrt{4-(d_{2})^2}\text{ cos}(\theta_{2}), \frac{d}{4d_{1}}-\\&\frac{(d_{1})^2+(d_{2})^2-(d_{3})^2}{4d_{1}d_{2}}\sqrt{4-(d_{2})^2}\text{ cos}(\theta_{2}), \frac{1}{2}\sqrt{4-(d_{2})^2}\text{ sin}(\theta_{2})\Big),\\
v_{5} =& \Big(\frac{(d_{1})^2-(d_{2})^2+(d_{3})^2}{2d_{1}}, \frac{d}{2d_{1}}, 0\Big),\\
v_{6} =& \Big(\frac{(d_{1})^2-(d_{2})^2+(d_{3})^2}{4d_{1}}-\frac{d}{4d_{1}d_{3}}\sqrt{4-(d_{3})^2}\text{ cos}(\theta_{3}), \frac{d}{4d_{1}}-\\&\frac{(d_{1})^2-(d_{2})^2+(d_{3})^2}{4d_{1}d_{3}}\sqrt{4-(d_{3})^2}\text{ cos}(\theta_{3}), \frac{1}{2}\sqrt{4-(d_{3})^2}\text{ sin}(\theta_{3})\Big),
\end{align*}

where $d=\sqrt{2(d_{1}d_{2})^2+2(d_{1}d_{3})^2+ 2(d_{2}d_{3})^2-(d_{1})^4-(d_{2})^4-(d_{3})^4}$.\\

Recall that the geometric knot invariant for hexagons, Joint Chirality-Curl, distinguishes between two types of both right-handed and left-handed trefoils with $curl(H)=\text{sign}((v_{3}-v_{1})\times(v_{5}-v_{1})\cdot(v_{2}-v_{1}))$. The following two lemmas give a relation between the curl of a hexagon and the possible dihedral angles.

\begin{lemma}\label{curl1}
	Let $H\in Equ_0(6)$. Let $H$ be parametrized using action-angle coordinates from the $T_{135}$ triangulation. If $H$ has Joint Chirality-Curl $(\pm1,1)$, then $\theta_i\in (0,\pi)$ for $i=1,2,3$.
\end{lemma}
\startproof 
Let $H\in Equ_{0}(6)$ be parametrized with action-angle coordinates $(d_1,d_2,d_3,\theta_1,\theta_2,\theta_3)$ arising from the $T_{135}$ triangulation. Let $H$ be in standard position. Since $v_1$, $v_3$, and $v_5$ are on the $xy$-plane oriented in a counter-clockwise direction, $curl(H)$ denotes the sign of the $z$-coordinate of $v_2$. Therefore if $curl(H)=1$, then $\theta_1\in(0,\pi)$. Suppose that $\theta_2, \theta_3\in (\pi, 2\pi)$. Then both  $e_4$ and $e_5$ lie below the $xy$-plane and can not pierce $T_2$. Thus $\Delta_2=0$ and $H$ has Joint Chirality-Curl $(0,0)$. If $\theta_2\in(\pi, 2\pi)$ and $\theta_3\in(0,\pi)$, then neither $e_6$ nor $e_1$ can pierce $T_4$ and $\Delta_4=0$. Similarly if $\theta_3\in(\pi, 2\pi)$ and $\theta_2\in(0,\pi)$, then $\Delta_6=0$. Therefore if $curl(H)=1$, then $\theta_{i}\in (0,\pi)$ for all $i\in{1,2,3}$.
\finishproof

\begin{lemma}\label{curl-1}
	Let $H\in Equ_0(6)$. Let $H$ be parametrized using action-angle coordinates from the $T_{135}$ triangulation. If $H$ has Joint Chirality-Curl $(\pm1,-1)$, then $\theta_i\in (\pi,2\pi)$ for $i=1,2,3$.
\end{lemma}
\startproof 
Let $H\in Equ_{0}(6)$ be parametrized with action-angle coordinates $(d_1,d_2,d_3,\theta_1,\theta_2,\theta_3)$ arising from the $T_{135}$ triangulation. Let $H$ be in standard position. If $curl(H)=-1$, then $\theta_1\in(\pi,2\pi)$. Similar to the previous argument of Lemma \ref{curl1}, if either of both of $\theta_2$ and $\theta_3$ are between $0$ and $\pi$, then $J(H)=(0,0)$. Therefore if $curl(H)=-1$ and $H$ is a trefoil, then $\theta_{i}\in (\pi, 2\pi)$ for all $i\in{1,2,3}$.
\finishproof\\


\subsection{Equilateral, Right-Handed, Positive Curl, Hexagonal Trefoils}

In this section, we will determine constraints on the values of $d_1, d_2, d_3, \theta_1, \theta_2,$ and $\theta_3$ in order to have a right-handed hexagonal trefoil with positive curl. 
First we will define a set of inequalities that must be satisfied in order for $H\in Equ_{0}(6)$ in standard position to have Joint Chirality-Curl $(1,1)$. 

\begin{prop}\label{disc}
	Let $H\in Equ_0(6)$ and parametrize $H$ with action-angle coordinates arising from the $T_{135}$ triangulation. If $J(H)=(1,1)$, then the following nine functions must be positive:

\begin{align*}
f_1= & d_2\sqrt{4-(d_2)^2}\text{sin}(\theta_2)\big(d_3d-((d_1)^2-(d_2)^2+(d_3)^2)\sqrt{4-(d_3)^2}\text{cos}(\theta_3)\big)-\\& d_3\sqrt{4-(d_3)^2}\text{sin}(\theta_3)\big(d_2d-((d_1)^2+(d_2)^2-(d_3)^2)\sqrt{4-(d_2)^2}\text{cos}(\theta_2)\big),\\
g_1= & \sqrt{4-(d_2)^2}\Big(\frac{-(d_1)^2+(d_2)^2+(d_3)^2}{2d_2 d_3}\text{cos}(\theta_2)\text{sin}(\theta_3)+\text{sin}(\theta_2)\text{cos}(\theta_3)\Big)\\ &-\frac{d\text{sin}(\theta_3)}{2d_3},\\
h_1= & \sqrt{4-(d_3)^2}\Big(\frac{-(d_1)^2+(d_2)^2+(d_3)^2}{2d_2 d_3}\text{cos}(\theta_3)\text{sin}(\theta_2)+\text{sin}(\theta_3)\text{cos}(\theta_2)\Big) \\ &-\frac{d\text{sin}(\theta_2)}{2d_2},\\
f_2 = & d_3\sqrt{4-(d_3)^2}\text{sin}(\theta_3)\big(d_1d-((d_1)^2+(d_2)^2-(d_3)^2)\sqrt{4-(d_1)^2}\text{cos}(\theta_1))- \\& d_1\sqrt{4-(d_1)^2}\text{sin}(\theta_1)\big(d_3d-(-(d_1)^2+(d_2)^2+(d_3)^2)\sqrt{4-(d_3)^2}\text{cos}(\theta_3)\big),\\
g_2 = &\sqrt{4-(d_3)^2}\Big(\frac{(d_1)^2-(d_2)^2+(d_3)^2}{2d_1 d_3}\text{cos}(\theta_3)\text{sin}(\theta_1)+\text{sin}(\theta_3)\text{cos}(\theta_1)\Big) \\ &-\frac{d\text{sin}(\theta_1)}{2d_1},\\
h_2 = &\sqrt{4-(d_1)^2}\Big(\frac{(d_1)^2-(d_2)^2+(d_3)^2}{2d_1 d_3}\text{cos}(\theta_1)\text{sin}(\theta_3)+\text{sin}(\theta_1)\text{cos}(\theta_3)\Big) \\ &-\frac{d\text{sin}(\theta_3)}{2d_3},\\
f_3 = & d_1\sqrt{4-(d_1)^2}\text{sin}(\theta_1)\big(d_2d-(-(d_1)^2+(d_2)^2+(d_3)^2)\sqrt{4-(d_2)^2}\text{cos}(\theta_2)\big)-\\&d_2\sqrt{4-(d_2)^2}\text{sin}(\theta_2)\big(d_1d-((d_1)^2-(d_2)^2+(d_3)^2)\sqrt{4-(d_1)^2}\text{cos}(\theta_1)\big),\\
g_3 = &\sqrt{4-(d_1)^2}\Big(\frac{(d_1)^2+(d_2)^2-(d_3)^2}{2d_1 d_2}\text{cos}(\theta_1)\text{sin}(\theta_2)+\text{sin}(\theta_1)\text{cos}(\theta_2)\Big)\\ & -\frac{d\text{sin}(\theta_2)}{2d_2},\\
h_3 = &\sqrt{4-(d_2)^2}\Big( \frac{(d_1)^2+(d_2)^2-(d_3)^2}{2d_1 d_2}\text{cos}(\theta_2)\text{sin}(\theta_1)+\text{sin}(\theta_2)\text{cos}(\theta_1)\Big) \\ & -\frac{d\text{sin}(\theta_1)}{2d_1}.
\end{align*}
\end{prop}

\startproof
Let $H\in Equ_0(6)$ and parametrize $H$ with action-angle coordinates from the $T_{135}$ triangulation. If $curl(H)=1$, then $\theta_i\in (0,\pi)$ for all $i$ by Lemma \ref{curl1}. Recall that for a hexagon, $H\in Equ_0(6)$ to be a right-handed trefoil then the algebraic intersection numbers $\Delta_i$ have to be equal to one for $i=2,4,6$. First we will consider $\Delta_4=1$. This means that the triangular disk $T_{4}$ containing $v_{3}, v_{4},$ and $v_{5}$ must be pierced by either the edge $e_{6}$ or $e_1$ so that the orientation on the edge agrees with the orientation on $T_4$ coming from a right-hand rule. If $\theta_2\in (0,\pi)$, then $e_6$ must pierce $T_4$ for $\Delta_4=1$. 

In order for the line going through $v_{1}$ and $v_{6}$ to pierce $T_{4}$ the following must be positive: 

\begin{align}
(v_{6}-v_{1})\times(v_{4}-v_{1})\centerdot(v_{3}-v_{1})&>0,\\
(v_{6}-v_{1})\times(v_{5}-v_{1})\centerdot(v_{4}-v_{1})&>0,\\
(v_{6}-v_{1})\times(v_{3}-v_{1})\centerdot(v_{5}-v_{1})&>0.
\end{align}

 \begin{figure}[h]
\centerline{\includegraphics[width=.5\textwidth]{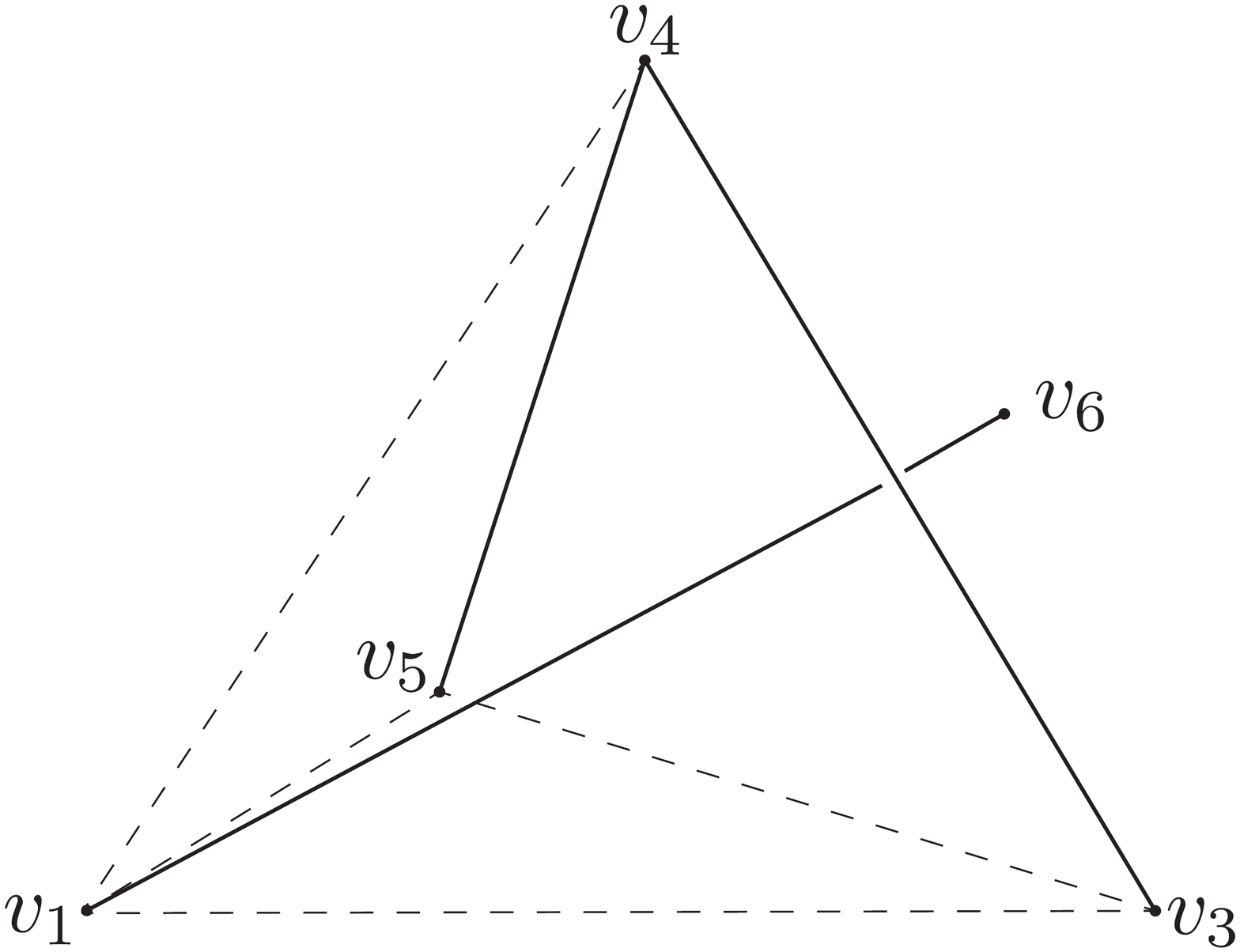}}
\caption{This figure shows the case where $e_6$ pierces $T_6$.}
\end{figure}

In addition, the plane containing $v_{3}, v_{4},$ and $v_{5}$ must separate $v_{1}$ from $v_{6}$. So the following must be negative 
\begin{align}
\big((v_{6}-v_{3})\times(v_{5}-v_{3})\centerdot(v_{4}-v_{3})\big)\big((v_{1}-v_{3})\times(v_{5}-v_{3})\centerdot(v_{4}-v_{3})\big)<0.
\end{align} 

Since $\theta_i\in(0,\pi)$ for all $i$, $(v_{6}-v_{1})\times(v_{3}-v_{1})\centerdot(v_{5}-v_{1})$ is always positive. Negating $(4)$ leaves three inequalities that must be satisfied so that $e_{6}$ pierces $T_4$. Letting $H$ be in standard position and evaluating the remaining three with action-angle coordinates gives three functions $f_1, g_1, h_1$ that must be positive for $H$ to have Joint Chirality-Curl $(1,1)$. Similarly, for $\Delta_2=1$ and $\Delta_6=1$ then $f_2, g_2, h_2$ and $f_3, g_3, h_3$ must be positive, respectively. Therefore if $H$ is in standard position and parametrized with action-angle coordinates form the $T_{135}$ triangulation, $f_i, g_i$, and $h_i$ must be positive for $H$ to have Joint Chirality-Curl $(1,1)$.
\finishproof\\

Next we will prove constraints on action-angle coordinates to get a right-handed, positive curl trefoil.

\begin{lemma}\label{equilateral}
Let $H\in Equ_{0}(6)$ and parametrize $H$ with with action-angle coordinates coming from the $T_{135}$ triangulation. If $H$ has Joint Chirality-Curl $(1,1)$, then the lengths of diagonals $d_i$ must be distinct. 
\end{lemma}
\startproof
Let $H$ be in $Equ_{0}(6)$. Parametrize $H$ with action-angle coordinates coming from the $T_{135}$ triangulation, so $H=(d_1,d_2,d_3,\theta_1,\theta_2,\theta_3)$ and $H$ is in standard position. Suppose $d_1=d_2=d_3=x$, for some $x\in(0,2)$. Then $v_1=(0,0,0), v_3=(x,0,0)$, and $v_5=(\frac{x}{2}, \frac{\sqrt{3}x}{2},0)$. First we will consider the case when $x=\sqrt{3}$. When $\theta_i=0$ for all $i$, $H$ is planar but singular with vertices $v_2$, $v_4$, and $v_6$ coinciding in a single point, $(\frac{\sqrt{3}}{2},\frac{1}{2},0)$. Additionally $e_1=e_6$, $e_2=e_3$, and $e_4=e_5$. As $\theta_1$ increases from $0$ to $\pi$, $v_{2}$ traverses a circle, $c_2$, of radius $\frac{1}{2}$ centered at $(\frac{\sqrt{3}}{2},0,0)$ lying in a plane parallel to the $yz$-plane. Similarly, as $\theta_3$ increases from $0$ to $2\pi$, $v_6$ moves along a circle, $c_6$, of radius $\frac{1}{2}$ centered at the midpoint of the edge connecting $v_1$ and $v_5$. Therefore $e_1$ sweeps out a circular cone, $C_{12}$, with vertex the origin and base circle $c_2$ and $e_6$ sweeps out a circular cone, $C_{61}$, with vertex the origin and base circle $c_6$. The circles $c_2$ and $c_6$ only intersect when $\theta_1=\theta_3=0$. Therefore the respective cones only intersect in the segment from $(0,0,0)$ to $(\frac{\sqrt{3}}{2},0,0)$, corresponding to edges $e_1$ and $e_6$ coinciding when $\theta_1=\theta_3=0$. This implies that $e_2$ can not pierce $T_6$. Thus if $x=\sqrt{3}$, $H$ can not have Joint Chirality-Curl $(1,1)$.

Next we consider the case when $\sqrt{3}<x<2$. When $\theta_1, \theta_2, \theta_3=0$, $H$ is planar and embedded. Hence $H$ is unknotted. Cones $C_{12}$ and $C_{61}$, formed by edges $e_1$ and $e_6$ as $\theta_1$ and $\theta_3$ vary, do not intersect. Therefore neither $T_2$ nor $T_6$ will be pierced by $H$, so $H$ will remain unknotted.

Now consider the case when $0<x<\sqrt{3}$. If $\theta_i=\text{cos}^{-1}(\frac{\sqrt{3}x}{3\sqrt{4-x^2}})$ for all $i$, then $v_2, v_4$, and $v_6$ coincide. If $\theta_i\in(\text{cos}^{-1}(\frac{\sqrt{3}x}{3\sqrt{4-x^2}}),\pi)$ for any $i\in \{1,2,3\}$ then $H$ will be unknotted. Therefore suppose that $\theta_i\in (0,\text{cos}^{-1}(\frac{\sqrt{3}x}{3\sqrt{4-x^2}}))$ for all $i$. If $\theta_1=\theta_3=0$, then $e_2$ and $e_5$ intersect in a point on the $xy$-plane. If $1\le x<\sqrt{3}$, then the point of intersection is interior of the triangle with vertices $v_1$, $v_3$, and $v_5$. As $\theta_1$ and $\theta_3$ increase from $0$ to $\text{cos}^{-1}(\frac{\sqrt{3}x}{3\sqrt{4-x^2}})$, $e_2$ and $e_5$ continue to intersect in a point. In order for $e_2$ to pierce $T_6$, then $\theta_3>\theta_1$. Similarly $e_1$ and $e_4$ intersect when $\theta_1=\theta_2$. In order for $e_4$ to pierce $T_2$, then $\theta_1>\theta_2$. When $\theta_2=\theta_3$, $e_6$ and $e_3$ intersect. In order for $e_6$ to pierce $T_4$, then $\theta_2>\theta_3$. This implies that $\theta_3>\theta_1>\theta_2>\theta_3$, a contradiction. When $0<x<1$ and $\theta_1=\theta_3=0$, then $e_2$ and $e_5$ intersect in a point that is exterior of the triangle with vertices $v_1$, $v_3$, and $v_5$. In order for $e_2$ to pierce $T_6$, then $\theta_3<\theta_1$. In order for $e_4$ to pierce $T_2$, then $\theta_1<\theta_2$. In order for $e_6$ to pierce $T_4$, then $\theta_2<\theta_3$. This implies that $\theta_3<\theta_1<\theta_2<\theta_3$, a contradiction. Therefore when $d_1=d_2=d_3$, $H$ can not have Joint Chirality-Curl $(1,1)$.
\finishproof\\

From Lemma \ref{curl1}, we know that all three dihedral angles must be in the interval $(0,\pi)$ for $curl(H)=1$. Next given any admissible triple of diagonal lengths, we prove a tighter constraints on the dihedral angles so that $H$ has Joint Chirality-Curl $(1,1)$.  

\begin{lemma}\label{angles}
Let $H\in Equ_{0}(6)$ and parametrize $H$ using action-angle coordinates with the $T_{135}$ triangulation. If $H$ has Joint Chirality-Curl $(1,1)$, then $\theta_{i}\in (0,\pi)$ for all $i\in{1,2,3}$ and $\theta_{1}+\theta_{2}<\pi$, $\theta_{1}+\theta_{3}<\pi$, and $\theta_{2}+\theta_{3}<\pi$. 
\end{lemma}
\startproof
Let $H\in Equ_{0}(6)$ be parametrized with action-angle coordinates\\ $(d_1,d_2,d_3,\theta_1,\theta_2,\theta_3)$ arising from the $T_{135}$ triangulation. If $H$ is a right-handed trefoil with $curl(H)=1$, then from Lemma \ref{curl1} we know $\theta_{i}\in (0,\pi)$ for all $i\in{1,2,3}$. If $\theta_i\in(0,\frac{\pi}{2})$ for all $i\in\{1,2,3\}$, then clearly $\theta_{1}+\theta_{2}<\pi$, $\theta_{1}+\theta_{3}<\pi$, and $\theta_{2}+\theta_{3}<\pi$. Additionally, if $\theta_i, \theta_j\in(0,\frac{\pi}{2})$, for any two distinct $i,j\in\{1,2,3\}$, then $\theta_i+\theta_j< \pi$.

Next we will show that if $\theta_1\in(\frac{\pi}{2},\pi)$, $\theta_2, \theta_3\in(0,\frac{\pi}{2})$ and $H$ has Joint Chirality-Curl $(1,1)$, then $\theta_{1}+\theta_{2}<\pi$ and $\theta_{1}+\theta_{3}<\pi$. Towards a contradiction, suppose that $\theta_1\in(\frac{\pi}{2},\pi)$, $\theta_2, \theta_3\in(0,\frac{\pi}{2})$, and $\theta_1+\theta_2= \pi$. Substituting $\theta_1=\pi-\theta_2$ into equation $g_3$ from Proposition \ref{disc} and using the facts that $\text{cos}(\pi-\theta_2)=-\text{cos}(\theta_2)$ and $\text{sin}(\pi-\theta_2)=\text{sin}(\theta_2)$, we obtain the following

\begin{align*}
g_3=\sqrt{4-(d_1)^2}\Big(-\frac{(d_1)^2+(d_2)^2-(d_3)^2}{2d_1 d_2}+1\Big)\text{sin}(\theta_2)\text{cos}(\theta_2)-\frac{d\text{sin}(\theta_2)}{2d_2}.
\end{align*}

We make the same substitutions into equation $h_3$ to obtain
\begin{align*}
h_3=\sqrt{4-(d_2)^2}\Big(\frac{(d_1)^2+(d_2)^2-(d_3)^2}{2d_1 d_2}-1\Big)\text{sin}(\theta_2)\text{cos}(\theta_2)-\frac{d\text{sin}(\theta_1)}{2d_1}.
\end{align*}

If $\frac{(d_1)^2+(d_2)^2-(d_3)^2}{2d_1 d_2}-1\ge0$, then $g_3$ is negative. Fix $d_1,d_2, d_3,$ and $\theta_2$, and now consider $g_3$ as a function of $\theta_1$. Since $\frac{(d_1)^2+(d_2)^2-(d_3)^2}{2d_1 d_2}-1\ge0$, then $\frac{(d_1)^2+(d_2)^2-(d_3)^2}{2d_1 d_2}>0$. This implies that the derivative of $g_3$ with respect to $\theta_1$, $$\sqrt{4-(d_1)^2}\Big(\frac{(d_1)^2+(d_2)^2-(d_3)^2}{2d_1 d_2}(-\text{sin}(\theta_1))\text{sin}(\theta_2)+\text{cos}(\theta_1)\text{cos}(\theta_2)\Big),$$ is negative for $\theta_1\in(\frac{\pi}{2},\pi)$ and $\theta_2\in(0,\frac{\pi}{2})$. Since $g_3$ is negative for $\theta_1=\pi-\theta_2$ and $g_3$ is decreasing, $g_3$ is negative for all $\theta_1>\pi-\theta_2$.
Next suppose $\frac{(d_1)^2+(d_2)^2-(d_3)^2}{2d_1 d_2}-1<0$, then $h_3$ is negative. This means that the plane, $P_2$, containing vertices $v_1, v_2$, and $v_3$ does not separate $v_4$ and $v_5$ when $\theta_1=\pi-\theta_2$. Therefore as $\theta_1$ increases to $\pi$ so that $\theta_1>\pi-\theta_2$, $P_2$ will not separate $v_4$ and $v_5$. Thus $h_3$ is negative for $\theta_1>\pi-\theta_2$.
Since both $g_3$ and $h_3$ must be positive for $H$ to a right-handed trefoil with $curl(H)=1$, we have reached a contradiction. Therefore if $\theta_1\in(0,\pi)$, $\theta_2, \theta_3\in(0,\frac{\pi}{2})$, and $\theta_1+\theta_2\ge \pi$, $H$ can not have Joint Chirality-Curl $(1,1)$.

Now suppose that $\theta_1\in(\frac{\pi}{2},\pi)$, $\theta_2, \theta_3\in(0,\frac{\pi}{2})$, and $\theta_1+\theta_3= \pi$. Similar to the previous argument, we will substitute $\theta_1=\pi-\theta_3$ into equation $g_2$ and use the facts that $\text{cos}(\pi-\theta_3)=-\text{cos}(\theta_3)$ and $\text{sin}(\pi-\theta_3)=\text{sin}(\theta_3)$. This results in the following:
\begin{align*}
g_2&=\sqrt{4-(d_3)^2}\Big(\frac{(d_1)^2-(d_2)^2+(d_3)^2}{2d_1 d_3} -1\Big)\text{cos}(\theta_3)\text{sin}(\theta_3)-\frac{d\text{sin}(\theta_1)}{2d_1}.
\end{align*}

Making the same substitutions into $h_2$ gives
\begin{align*}
h_2&=\sqrt{4-(d_1)^2}\Big(-\frac{(d_1)^2-(d_2)^2+(d_3)^2}{2d_1 d_3} +1\Big)\text{cos}(\theta_3)\text{sin}(\theta_3)-\frac{d\text{sin}(\theta_3)}{2d_3}.
\end{align*}

If $\frac{(d_1)^2-(d_2)^2+(d_3)^2}{2d_1 d_3} -1\le 0$, then $g_2$ is negative. If $\frac{(d_1)^2-(d_2)^2+(d_3)^2}{2d_1 d_3} -1> 0$, then $h_2$ is negative. Since both equations must be positive to have Joint Chirality-Curl $(1,1)$, we have reached a contradiction. Therefore if $H$ is a right-handed trefoil with $curl(H)=1$ and $\theta_1\in(\frac{\pi}{2},\pi)$, $\theta_2, \theta_3\in(0,\frac{\pi}{2})$, then $\theta_{1}+\theta_{2}<\pi$ and $\theta_{1}+\theta_{3}<\pi$.

The cases when $\theta_2\in(\frac{\pi}{2}, \pi)$, $\theta_1, \theta_3\in(0,\frac{\pi}{2})$ and $\theta_3\in(\frac{\pi}{2}, \pi)$, $\theta_1, \theta_2\in(0,\frac{\pi}{2})$ are proven in the same manner. Thus if $H$ has Joint Chirality-Curl $(1,1)$, then $\theta_{i}\in (0,\pi)$ for all $i\in{1,2,3}$ and $\theta_{1}+\theta_{2}<\pi$, $\theta_{1}+\theta_{3}<\pi$, and $\theta_{2}+\theta_{3}<\pi$.

\finishproof\\

For $H\in Equ_{0}(6)$ to be a right-handed trefoil with $curl(H)=1$, only one of the dihedral angles can be greater than $\pi/2$ with the additional condition that the sum of any two angles must be less than $\pi$. This portion of the cube $[0,2\pi]^3$ is shown in the following Figure \ref{torusportion}.

\begin{figure}[h]
\centerline{\includegraphics[width=.7\textwidth]{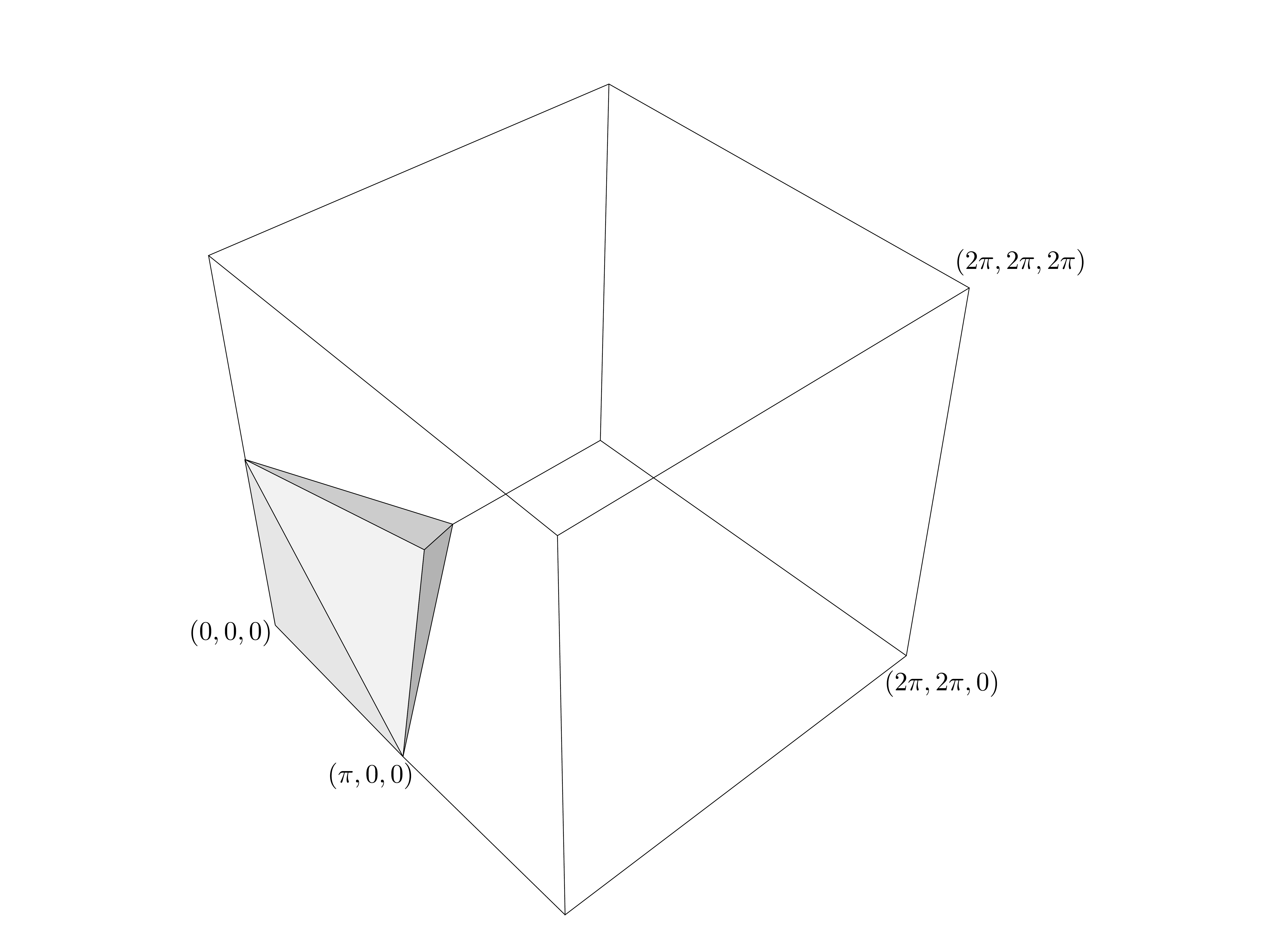}}
\caption{The corresponding angles from Lemma \ref{angles} for an equilateral hexagon to have Joint Chirality-Curl $(1,1)$.}
\label{torusportion}
\end{figure}

When all three diagonals from the $T_{135}$ triangulation have equal length, $H$ can not have Joint Chirality-Curl $(1,1)$. So, we continue our analysis with the case where two of the diagonals have equal lengths.\\

 coincide when $\theta_1=\theta_2=0$. This occurs when $d_3=d_1\sqrt{4-(d_1)^2}$. If $d_3\ge d_1\sqrt{4-(d_1)^2}$, then $H$ has Joint Chirality $(0,0)$ for all $\theta_i$. Now we will consider the different ranges of $\theta_i$ for $H$ to have Joint Chirality-Curl $(1,1)$ from Lemma \ref{angles}.

Next we consider the case where the three diagonals of the $T_{135}$ triangulation are distinct.\\ 
\begin{lemma}\label{polytopeportionR+}
Let $H\in Equ_{0}(6)$ and parametrize $H$ using action-angle coordinates with the $T_{135}$ triangulation. Suppose $d_{1}, d_{2},$ and $d_{3}$ are distinct and let $d_{i}>d_{j}, d_{k}$. If $J(H)=(1,1)$ then $\theta_{i}\in(0,\pi)$ and $\theta_{j},\theta_{k}\in(0,\pi/2)$. Moreover, if $d_{i}>\sqrt{(d_{j})^2+(d_{k})^2}$ then $\theta_{i}\in(\pi/2,\pi)$.\\
\end{lemma}
\startproof
Let $H\in Equ_{0}(6)$ be in standard position so that $v_1$, $v_3$, and $v_5$ are on the $xy$-plane. Suppose that the lengths of the diagonals are distinct and that $d_2>d_3>d_1$. We will show that if $H$ has Joint Chirality-Curl $(1,1)$ then $\theta_2\in (0,\pi)$ and $\theta_1,\theta_3\in(0,\frac{\pi}{2})$. Let $l_2$ be the line in the $xy$-plane perpendicular to the segment connecting $v_1$ and $v_3$ intersecting at the midpoint, $m_2$. Similarly, we define $l_4$ and $l_6$ for segments connecting $v_3$ to $v_5$ and $v_5$ to $v_1$ respectively. The three lines intersect in a unique point, $k$, the circumcenter of the triangle with vertices $v_1$, $v_3$, and $v_5$. Moreover, $l_i$ represents the orthogonal projection of $v_i$ onto the $xy$-plane as $\theta_i$ varies and $k$ is the projection of where all vertices coincide, if such point exists.

\begin{figure}[h]
\centerline{\includegraphics[width=.55\textwidth]{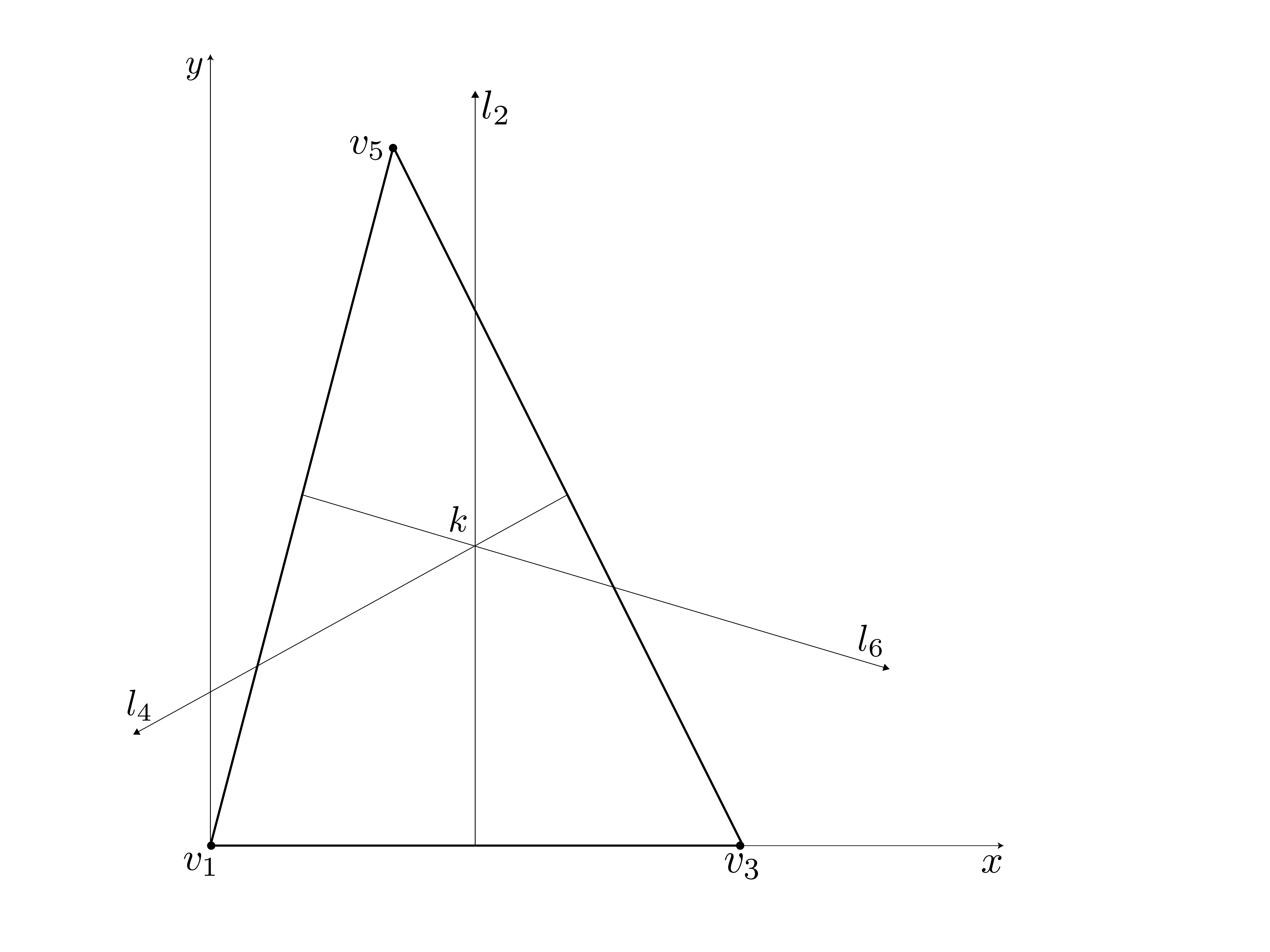}\includegraphics[width=.55\textwidth]{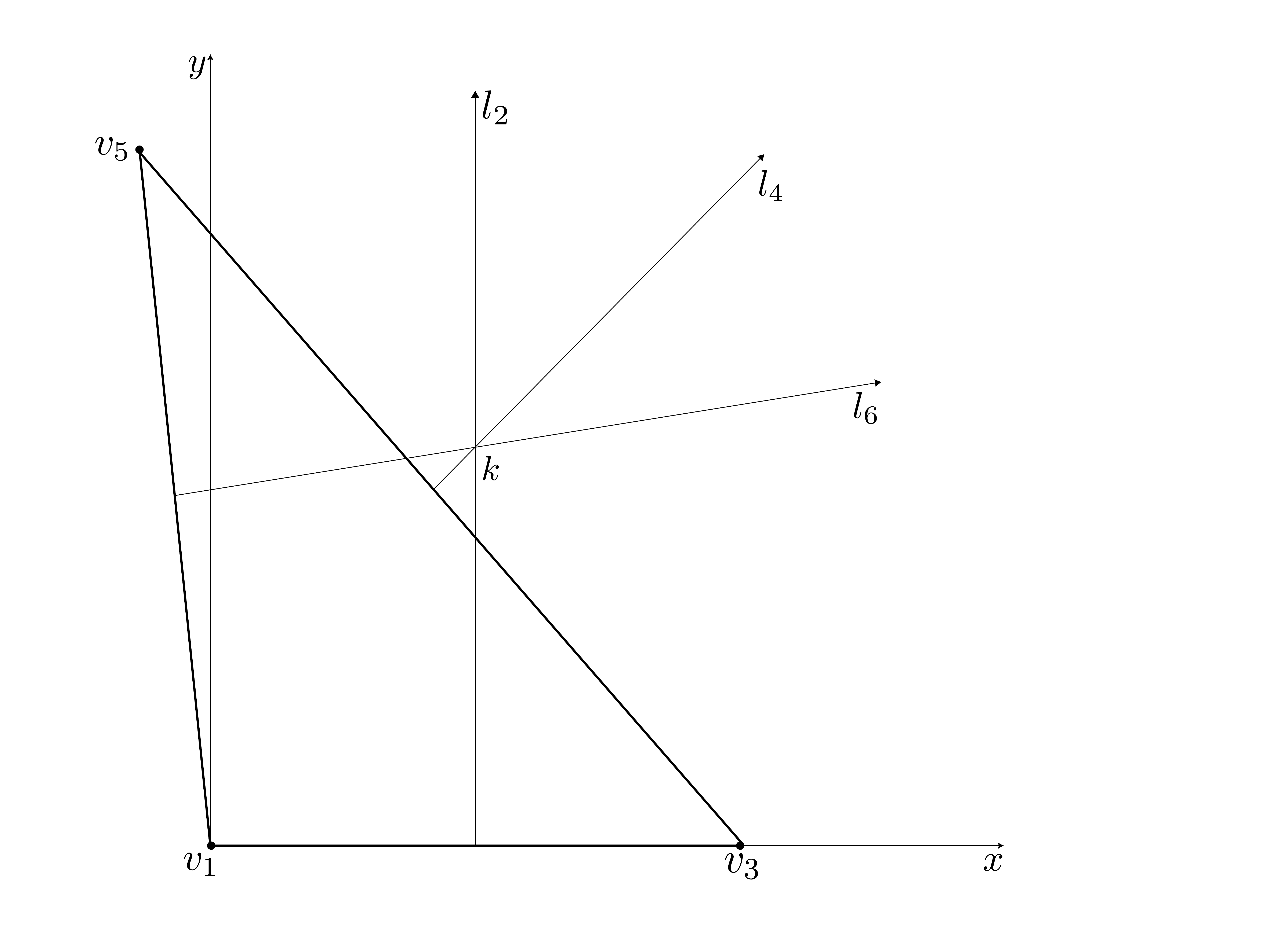}}
\caption{The figures show the triangle spanned by vertices $(v_1 ,v_3 ,v_5 )$ with perpendicular bisector $l_i$. In the figure on the left, $(d_2)^2 < (d_1)^2+(d_3)^2$ and so $k$ is interior of the triangle. In the figure on the right, $(d_2)^2 >(d_1)^2+(d_3)^2$ and so $k$ is exterior of the triangle.}
\label{obtuse}
\end{figure}

Since $d_3>d_1$ then $l_4$ intersects the segment connecting $v_1$ to $v_5$ instead of the segment connecting $v_1$ and $v_3$. Suppose towards contradiction that $\theta_1\in(\frac{\pi}{2},\pi)$. Then the plane perpendicular to the $xy$-plane containing $l_4$ separates $e_4$ and $T_2$. Therefore $H$ can not have Joint Chirality-Curl $(1,1)$ if $\theta_1\in(\frac{\pi}{2},\pi)$. Let $\phi_1$ be the angle for $\theta_1$ where $v_2$ projects onto $k$. If $\theta_1\in (\phi_1,\frac{\pi}{2})$ then $e_4$ and $T_2$ are still separated by the plane through $l_4$. Therefore $\theta_1\in (0,\phi_1)$. 
Next suppose that $\theta_3\in(\frac{\pi}{2},\pi)$. Since $d_2>d_3$ the $l_2$ intersects the segment connecting $v_3$ and $v_5$. This means the plane perpendicular to the $xy$-plane containing $l_2$ separates $e_2$ and $T_6$. Therefore we have reached a contradiction and $\theta_3\in (0,\frac{\pi}{2})$. Let $\phi_3$ be the angle for $\theta_3$ for which $v_6$ projects onto $k$. In order for $e_2$ to pierce $T_6$ then $\theta_3\in (0,\phi_3)$. Let $\phi_2$ be the angle for $\theta_2$ for which $v_4$ projects onto $k$. Let $p$ be the point where $e_1$ intersects $e_4$ when $\theta_1=0$ and $\theta_2=\pi$. In order for $e_4$ to intersect $T_2$, $e_4$ must intersect the cone spanned by $e_1$. The two cones will intersect along an arc connecting $p$ to point which projects onto $k$. If $\theta_2<\phi_2$, then $e_4$ no longer intersects the cone spanned by $e_1$. Then $\theta_2\in (\phi_2,\pi)$ for $H$ to have Joint Chirality-Curl $(1,1)$. If $d_2>\sqrt{(d_1)^2+(d_3)^2}$ then the triangle with vertices $v_1$, $v_3$, and $v_5$ is obtuse, as shown in Figure \ref{obtuse}. Therefore $k$ is exterior of the triangle and $\phi_2>\frac{\pi}{2}$. Hence if $d_2>\sqrt{(d_1)^2+(d_3)^2}$ then $\theta_1,\theta_3\in (0,\frac{\pi}{2})$ and $\theta_2\in (\frac{\pi}{2},\pi)$.
\finishproof\\

\begin{figure}[h]
\centerline{\includegraphics[width=.5\textwidth]{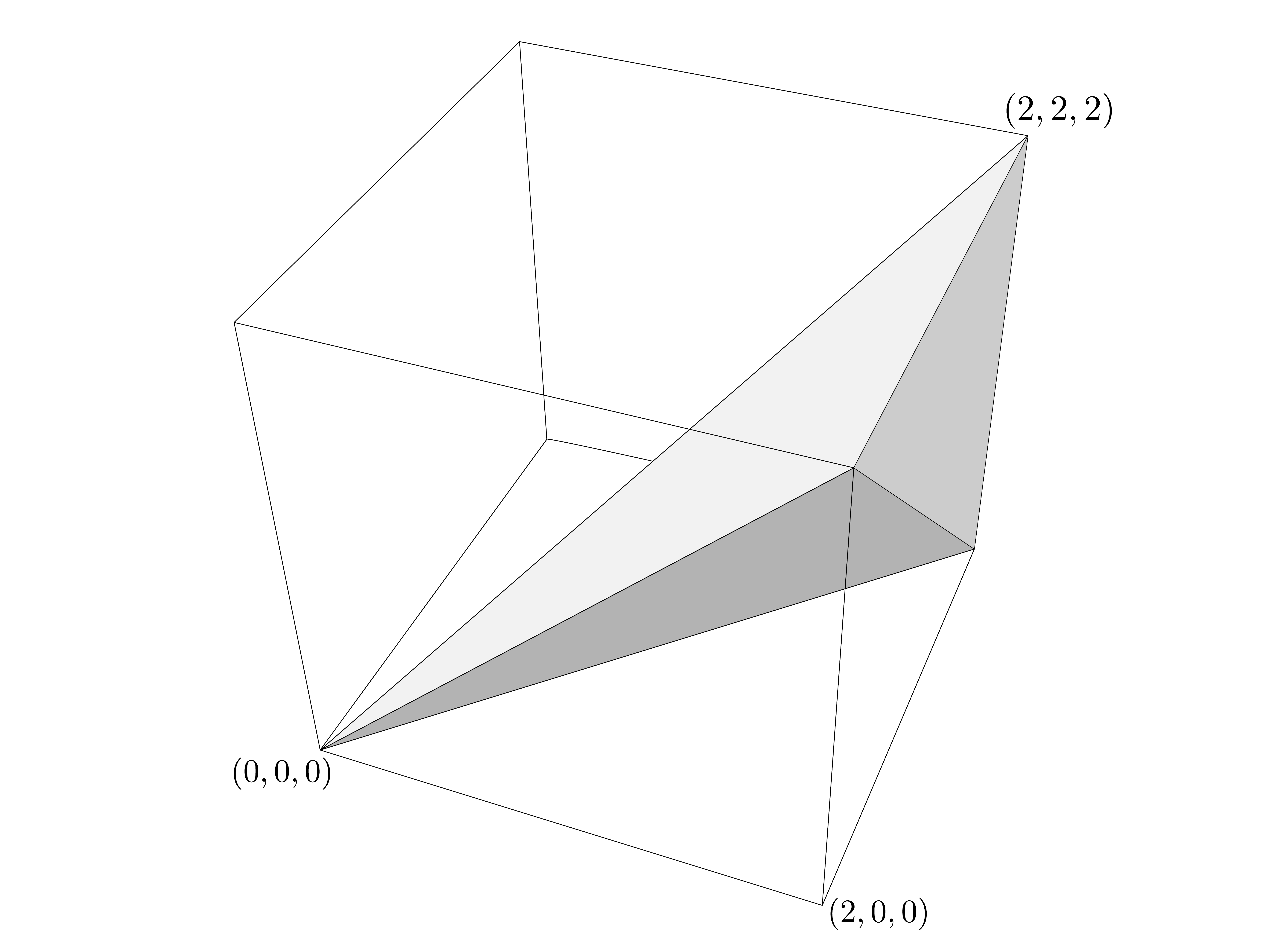}\includegraphics[width=.5\textwidth]{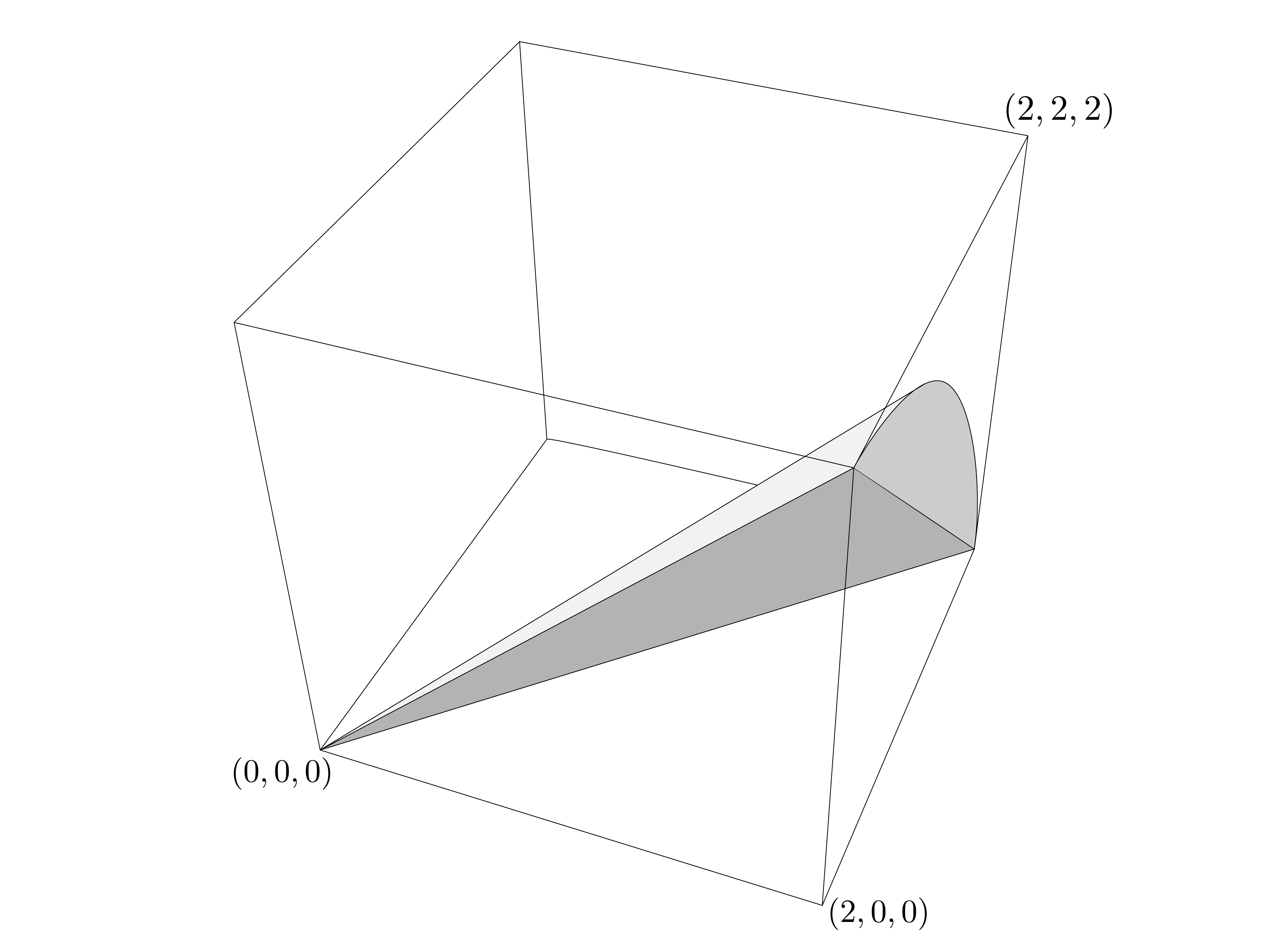}}
\caption{The figure on the left shows the portion of the moment polytope, $P_6$, where $d_1>d_2, d_3$. The figure on the left shows the portion of the moment polytope where additionally $(d_1)^2>(d_2)^2+(d_3)^2.$}
\label{fig:dividepolytope}
\end{figure}

The moment polytope corresponding to the $T_{135}$ triangulation is split into three equal regions, depending on the which diagonal length is largest. The function $d_1=\sqrt{(d_2)^2+(d_3)^2}$ divides the third of the polytope where $d_1>d_2, d_3$ into two regions, one for acute triangles and one for obtuse triangles, as shown in Figure \ref{fig:dividepolytope}.

\subsection{Equilateral, Left-Handed, Positive Curl, Hexagonal Trefoils}

In this section, we will discuss constraints for $H$ to be a left-handed hexagonal trefoil with positive curl. 

\begin{prop}\label{discL+}
Let $H\in Equ_0(6)$ and parametrize $H$ with action-angle coordinates arising from the $T_{135}$ triangulation. If $J(H)=(-1,1)$, then $f_i<0$, $g_i>0$, and $h_i>0$, for all $i$, where
	
\begin{align*}
f_1= & d_2\sqrt{4-(d_2)^2}\text{sin}(\theta_2)\big(d_3d-((d_1)^2-(d_2)^2+(d_3)^2)\sqrt{4-(d_3)^2}\text{cos}(\theta_3)\big)-\\& d_3\sqrt{4-(d_3)^2}\text{sin}(\theta_3)\big(d_2d-((d_1)^2+(d_2)^2-(d_3)^2)\sqrt{4-(d_2)^2}\text{cos}(\theta_2)\big),\\
g_1= & \sqrt{4-(d_2)^2}\Big(\frac{-(d_1)^2+(d_2)^2+(d_3)^2}{2d_2 d_3}\text{cos}(\theta_2)\text{sin}(\theta_3)+\text{sin}(\theta_2)\text{cos}(\theta_3)\Big)\\ &-\frac{d\text{sin}(\theta_3)}{2d_3},\\
h_1= & \sqrt{4-(d_3)^2}\Big(\frac{-(d_1)^2+(d_2)^2+(d_3)^2}{2d_2 d_3}\text{cos}(\theta_3)\text{sin}(\theta_2)+\text{sin}(\theta_3)\text{cos}(\theta_2)\Big) \\ &-\frac{d\text{sin}(\theta_2)}{2d_2},\\
f_2 = & d_3\sqrt{4-(d_3)^2}\text{sin}(\theta_3)\big(d_1d-((d_1)^2+(d_2)^2-(d_3)^2)\sqrt{4-(d_1)^2}\text{cos}(\theta_1))- \\& d_1\sqrt{4-(d_1)^2}\text{sin}(\theta_1)\big(d_3d-(-(d_1)^2+(d_2)^2+(d_3)^2)\sqrt{4-(d_3)^2}\text{cos}(\theta_3)\big),\\
g_2 = &\sqrt{4-(d_3)^2}\Big(\frac{(d_1)^2-(d_2)^2+(d_3)^2}{2d_1 d_3}\text{cos}(\theta_3)\text{sin}(\theta_1)+\text{sin}(\theta_3)\text{cos}(\theta_1)\Big) \\ &-\frac{d\text{sin}(\theta_1)}{2d_1},\\
h_2 = &\sqrt{4-(d_1)^2}\Big(\frac{(d_1)^2-(d_2)^2+(d_3)^2}{2d_1 d_3}\text{cos}(\theta_1)\text{sin}(\theta_3)+\text{sin}(\theta_1)\text{cos}(\theta_3)\Big) \\ &-\frac{d\text{sin}(\theta_3)}{2d_3},\\
f_3 = & d_1\sqrt{4-(d_1)^2}\text{sin}(\theta_1)\big(d_2d-(-(d_1)^2+(d_2)^2+(d_3)^2)\sqrt{4-(d_2)^2}\text{cos}(\theta_2)\big)-\\&d_2\sqrt{4-(d_2)^2}\text{sin}(\theta_2)\big(d_1d-((d_1)^2-(d_2)^2+(d_3)^2)\sqrt{4-(d_1)^2}\text{cos}(\theta_1)\big),\\
g_3 = &\sqrt{4-(d_1)^2}\Big(\frac{(d_1)^2+(d_2)^2-(d_3)^2}{2d_1 d_2}\text{cos}(\theta_1)\text{sin}(\theta_2)+\text{sin}(\theta_1)\text{cos}(\theta_2)\Big)\\ & -\frac{d\text{sin}(\theta_2)}{2d_2},\\
h_3 = &\sqrt{4-(d_2)^2}\Big( \frac{(d_1)^2+(d_2)^2-(d_3)^2}{2d_1 d_2}\text{cos}(\theta_2)\text{sin}(\theta_1)+\text{sin}(\theta_2)\text{cos}(\theta_1)\Big) \\ & -\frac{d\text{sin}(\theta_1)}{2d_1}.
\end{align*}	
\end{prop}
\startproof 
Let $H\in Equ_0(6)$ and parametrize $H$ using action-angle coordinates for the $T_{135}$ triangulation. If $curl(H)=1$, then $\theta_i\in (0,\pi)$ for all $i$. Additionally, if $H$ has Joint Chirality-Curl $(-1,1)$ the all algebraic intersection numbers $\Delta_i$ must be negative. First we will consider that condition that $\Delta_4=-1$, meaning the algebraic intersection of $T_4$ with $H$ is $-1$. If $\theta_2\in (0,\pi)$ this means that $e_1$ must pierce $T_4$. Therefore the line going through $v_1$ and $v_2$ must pass through $T_4$ and the following three inequalities must be satisfied:
\begin{align}
	(v_2-v_1)\times(v_4-v_1)\centerdot(v_3-v_1)>0,\\
	(v_2-v_1)\times(v_5-v_1)\centerdot(v_4-v_1)>0,\\
	(v_2-v_1)\times(v_3-v_1)\centerdot(v_5-v_1)>0.
\end{align}

Since $\theta_i\in(0,\pi)$, then $(v_2-v_1)\times(v_3-v_1)\centerdot(v_5-v_1)$ is always positive. Suppose $H$ is in standard position. Evaluating the remaining two expressions with action-angle coordinates from the $T_{135}$ triangulation results in $h_3>0$ and $f_3<0$.

Additionally, the plane containing $v_3, v_4$, and $v_5$ must separate $v_1$ and $v_2$. Therefore the following must be negative:
$$\big( (v_2-v_3)\times(v_5-v_3)\centerdot(v_4-v_3)\big)\centerdot\big( (v_1-v_3)\times(v_5-v_3)\centerdot (v_4-v_3)\big)<0.$$ 
This constraint is equivalent to $g_3>0$. 
 
Similarly, the conditions that the $\Delta_2=-1$ and $\Delta_6=-1$ are equivalent to $f_2<0, g_2>0, h_2>0$ and $f_3<0, g_3>0, h_3>0$, respectively. Therefore if $H\in Equ_0(6)$ is in standard position and has Joint Chirality-Curl $(-1,1)$, then $f_i<0$, $g_i>0$, and $h_i>0$ for all $i$.
\finishproof\\

Next we define possible dihedral angles for an equilateral hexagon to have Joint Chirality-Curl $(-1,1)$.

\begin{lemma}\label{angles(-1,1)}
Let $H\in Equ_{0}(6)$ and parametrize $H$ using action-angle coordinates with the $T_{135}$ triangulation. If $H$ has Joint Chirality-Curl $(-1,1)$, then $\theta_{i}\in (0,\pi)$ for all $i\in{1,2,3}$ and $\theta_{1}+\theta_{2}<\pi$, $\theta_{1}+\theta_{3}<\pi$, and $\theta_{2}+\theta_{3}<\pi$. 
\end{lemma}
\startproof
Let $H\in Equ_{0}(6)$ be parametrized with action-angle coordinates $(d_1,d_2,d_3,\theta_1,\theta_2,\theta_3)$ arising from the $T_{135}$ triangulation. If $H$ has Joint Chirality-Curl $(-1,1)$, then from Lemma \ref{curl1} we know $\theta_{i}\in (0,\pi)$ for all $i\in{1,2,3}$. If $\theta_i\in(0,\frac{\pi}{2})$ for all $i\in\{1,2,3\}$, then clearly $\theta_{1}+\theta_{2}<\pi$, $\theta_{1}+\theta_{3}<\pi$, and $\theta_{2}+\theta_{3}<\pi$. Additionally if $\theta_i, \theta_j\in(0,\frac{\pi}{2})$, for any two distinct $i,j\in\{1,2,3\}$, then $\theta_i+\theta_j< \pi$.

Next we will show that if $\theta_1\in(\frac{\pi}{2},\pi)$, $\theta_2, \theta_3\in(0,\frac{\pi}{2})$ and $H$ has Joint Chirality-Curl $(-1,1)$, then $\theta_{1}+\theta_{2}<\pi$ and $\theta_{1}+\theta_{3}<\pi$. Towards a contradiction, suppose that $\theta_1\in(\frac{\pi}{2},\pi)$, $\theta_2, \theta_3\in(0,\frac{\pi}{2})$, and $\theta_1+\theta_2= \pi$. Substituting $\theta_1=\pi-\theta_2$ into equation $g_3$ from Proposition \ref{discL+} and using the facts that $\text{cos}(\pi-\theta_2)=-\text{cos}(\theta_2)$ and $\text{sin}(\pi-\theta_2)=\text{sin}(\theta_2)$, we obtain the following
\begin{align*}
g_3&=\sqrt{4-(d_1)^2}\Big(-\frac{(d_1)^2+(d_2)^2-(d_3)^2}{2d_1 d_2}+1\Big)\text{sin}(\theta_2)\text{cos}(\theta_2)-\frac{d\text{sin}(\theta_2)}{2d_2}.
\end{align*}

We make the same substitutions into equation $h_3$ to obtain

\begin{align*}
h_3&=\sqrt{4-(d_2)^2}\Big(\frac{(d_1)^2+(d_2)^2-(d_3)^2}{2d_1 d_2}-1\Big)\text{sin}(\theta_2)\text{cos}(\theta_2)-\frac{d\text{sin}(\theta_1)}{2d_1}.
\end{align*}

If $\frac{(d_1)^2+(d_2)^2-(d_3)^2}{2d_1 d_2}-1\ge0$, then $g_3$ is negative. Fix $d_1,d_2, d_3,$ and $\theta_2$, and now consider $g_3$ as a function of $\theta_1$. Since $\frac{(d_1)^2+(d_2)^2-(d_3)^2}{2d_1 d_2}-1\ge0$, then $\frac{(d_1)^2+(d_2)^2-(d_3)^2}{2d_1 d_2}>0$. This implies that the derivative of $g_3$ with respect to $\theta_1$, $$\sqrt{4-(d_1)^2}\Big(\frac{(d_1)^2+(d_2)^2-(d_3)^2}{2d_1 d_2}(-\text{sin}(\theta_1))\text{sin}(\theta_2)+\text{cos}(\theta_1)\text{cos}(\theta_2)\Big),$$ is negative for $\theta_1\in(\frac{\pi}{2},\pi)$ and $\theta_2\in(0,\frac{\pi}{2})$. Since $g_3$ is negative for $\theta_1=\pi-\theta_2$ and $g_3$ is decreasing, $g_3$ is negative for all $\theta_1>\pi-\theta_2$.
Next suppose $\frac{(d_1)^2+(d_2)^2-(d_3)^2}{2d_1 d_2}-1<0$, then $h_3$ is negative. Thus $h_3$ is negative for $\theta_1>\pi-\theta_2$.
Since both $g_3$ and $h_3$ must be positive for $H$ to have Joint Chirality-Curl $(-1,1)$, we have reached a contradiction. Therefore if $\theta_1\in(0,\pi)$, $\theta_2, \theta_3\in(0,\frac{\pi}{2})$, and $\theta_1+\theta_2\ge \pi$, $H$ can not have Joint Chirality-Curl $(-1,1)$.

Now suppose that $\theta_1\in(\frac{\pi}{2},\pi)$, $\theta_2, \theta_3\in(0,\frac{\pi}{2})$, and $\theta_1+\theta_3= \pi$. Similar to the previous argument, we will substitute $\theta_1=\pi-\theta_3$ into equation $g_2$ and use the facts that $\text{cos}(\pi-\theta_3)=-\text{cos}(\theta_3)$ and $\text{sin}(\pi-\theta_3)=\text{sin}(\theta_3)$. This results in the following:
\begin{align*}
g_2&=\sqrt{4-(d_3)^2}\Big(\frac{(d_1)^2-(d_2)^2+(d_3)^2}{2d_1 d_3} -1\Big)\text{cos}(\theta_3)\text{sin}(\theta_3)-\frac{d\text{sin}(\theta_1)}{2d_1}.
\end{align*}

Making the same substitutions into $h_2$ gives
\begin{align*}
h_2&=\sqrt{4-(d_1)^2}\Big(-\frac{(d_1)^2-(d_2)^2+(d_3)^2}{2d_1 d_3} +1\Big)\text{cos}(\theta_3)\text{sin}(\theta_3)-\frac{d\text{sin}(\theta_3)}{2d_3}.
\end{align*}

If $\frac{(d_1)^2-(d_2)^2+(d_3)^2}{2d_1 d_3} -1\le 0$, then $g_2$ is negative. If $\frac{(d_1)^2-(d_2)^2+(d_3)^2}{2d_1 d_3} -1> 0$, then $h_2$ is negative. Since both equations must be positive to have Joint Chirality-Curl $(-1,1)$, we have reached a contradiction. Therefore if $H$ is a left-handed trefoil with $curl(H)=1$ and $\theta_1\in(\frac{\pi}{2},\pi)$, $\theta_2, \theta_3\in(0,\frac{\pi}{2})$, then $\theta_{1}+\theta_{2}<\pi$ and $\theta_{1}+\theta_{3}<\pi$. 

Again the cases when $\theta_2\in(\frac{\pi}{2}, \pi)$, $\theta_1, \theta_3\in(0,\frac{\pi}{2})$ and $\theta_3\in(\frac{\pi}{2}, \pi)$, $\theta_1, \theta_2\in(0,\frac{\pi}{2})$ are proven in the same manner. Thus if $H$ has Joint Chirality-Curl $(-1,1)$, then $\theta_{i}\in (0,\pi)$ for all $i\in{1,2,3}$ and $\theta_{1}+\theta_{2}<\pi$, $\theta_{1}+\theta_{3}<\pi$, and $\theta_{2}+\theta_{3}<\pi$.

\finishproof\\

Next we consider the case where the three diagonals of the $T_{135}$ triangulation are distinct.\\ 
\begin{lemma}\label{polytopeportionL+}
Let $H\in Equ_{0}(6)$ and parametrize $H$ using action-angle coordinates with the $T_{135}$ triangulation. Suppose $d_{1}, d_{2},$ and $d_{3}$ are distinct and let $d_{i}>d_{j}, d_{k}$. If $J(H)=(-1,1)$ then $\theta_{i}\in(0,\pi)$ and $\theta_{j},\theta_{k}\in(0,\pi/2)$. Moreover, if $d_{i}>\sqrt{(d_{j})^2+(d_{k})^2}$ then $\theta_{i}\in(\pi/2,\pi)$.\\
\end{lemma}
\startproof
Let $H\in Equ_{0}(6)$ be in standard position so that $v_1$, $v_3$, and $v_5$ are on the $xy$-plane. Suppose that the lengths of the diagonals are distinct and that $d_2>d_1>d_3$. We will show that if $H$ has Joint Chirality-Curl $(-1,1)$ then $\theta_2\in (0,\pi)$ and $\theta_1,\theta_3\in(0,\frac{\pi}{2})$. Let $l_2$ be perpendicular bisector to the segment connecting $v_1$ and $v_3$. Similarly, we define $l_4$ and $l_6$ to be the perpendicular bisectors to segments connecting $v_3$ to $v_5$ and $v_5$ to $v_1$, respectively. The three lines intersect in a unique point, $k$, the circumcenter of the triangle spanned by $(v_1, v_3, v_5)$. The orthogonal projection of $v_i$ onto the $xy$-plane lies on $l_i$. In addition, $k$ is the orthogonal projection of where all vertices coincide, if such point exists. 

Since $d_1>d_3$ then $l_4$ intersects the segment connecting $v_1$ to $v_3$ instead of the segment connecting $v_1$ and $v_5$. Suppose towards contradiction that $\theta_3\in(\frac{\pi}{2},\pi)$. Then the plane perpendicular to the $xy$-plane containing $l_4$ separates $e_3$ and $T_6$. Therefore $H$ can not have Joint Chirality-Curl $(-1,1)$ if $\theta_3\in(\frac{\pi}{2},\pi)$.

Next suppose that $\theta_1\in(\frac{\pi}{2},\pi)$. Since $d_2>d_1$ the $l_6$ intersects the segment connecting $v_3$ and $v_5$. This means the plane perpendicular to the $xy$-plane containing $l_6$ separates $e_5$ and $T_2$. Therefore we have reached a contradiction and $\theta_1\in (0,\frac{\pi}{2})$.

 Let $\psi_3$ be the angle for $\theta_3$ for which $v_6$ projects onto $k$. If $\theta_3\in(\psi_3,\pi)$, then the plane perpendicular to the $xy$-plane containing $l_4$ still separates $e_3$ and $T_6$. Thus if $H$ has Joint Chirality-Curl $(-1,1)$, then $\theta_3\in(0,\psi_3)$. Let $\psi_1$ be the angle of $\theta_1$ so that $v_1$ projects onto $k$. If $e_5$ is to intersect $T_2$, then $\theta_1\in(0,\psi_1)$. Let $\psi_2$ be the angle for $\theta_2$ for which $v_4$ projects onto $k$. Since $\theta_1\in (0,\psi_1)$ and $\theta_3\in (0,\psi_3)$ then $\theta_2\in (\psi_2,\pi)$ for $H$ to have Joint Chirality-Curl $(-1,1)$. If $d_2>\sqrt{(d_1)^2+(d_3)^2}$ then the triangle spanned by $(v_1, v_3, v_5)$ is obtuse. Therefore $k$ is exterior of the triangle spanned by $(v_1, v_3, v_5)$ and $\psi_2>\frac{\pi}{2}$. Hence if $d_2>\sqrt{(d_1)^2+(d_3)^2}$ then $\theta_1,\theta_3\in (0,\frac{\pi}{2})$ and $\theta_2\in (\frac{\pi}{2},\pi)$.
\finishproof\\

\section{Knotting Probability of Hexagonal Trefoils} 

In this section, we will discuss the probability that a random equilateral hexagon is knotted. It has been proven that at least $\frac{1}{3}$ of hexagons with total length $2$ are unknotted\cite{JC}. Using action-angle coordinates and Calvo's geometric invariant $curl$, Cantarella and Shonkwiler \cite{JC} prove that at least $\frac{1}{2}$ of the space of equilateral hexagons consists of unknots. In order to gain intuition on the tightness of these bounds, we performed a Monte Carlo experiment. We randomly sampled a point in the moment polytope $P_6$ and a point in the cube $[0,2\pi]^3$. We then tested whether this point satisfies the necessary constraints to have a hexagonal trefoil with Joint Chirality-Curl $(1,1)$, described in Proposition \ref{disc}. Taking a sample size of $10$ million configurations, repeating this experiment multiple times, we found that on average the fraction of $(1,1)$ trefoils is $3.426005\times 10^{-5}$ with standard deviation $2.241511 \times 10^{-6}$. Since there are four types of trefoils, we estimate that the knotting probability for equilateral hexagons is about $1.370402 \times 10^{-4}$. Using the lemmas from the previous section, we improve the theoretical bound.

\begin{figure}[h]
\centerline{\includegraphics[width=.5\textwidth]{polytopeobtuse.pdf}\includegraphics[width=.45\textwidth]{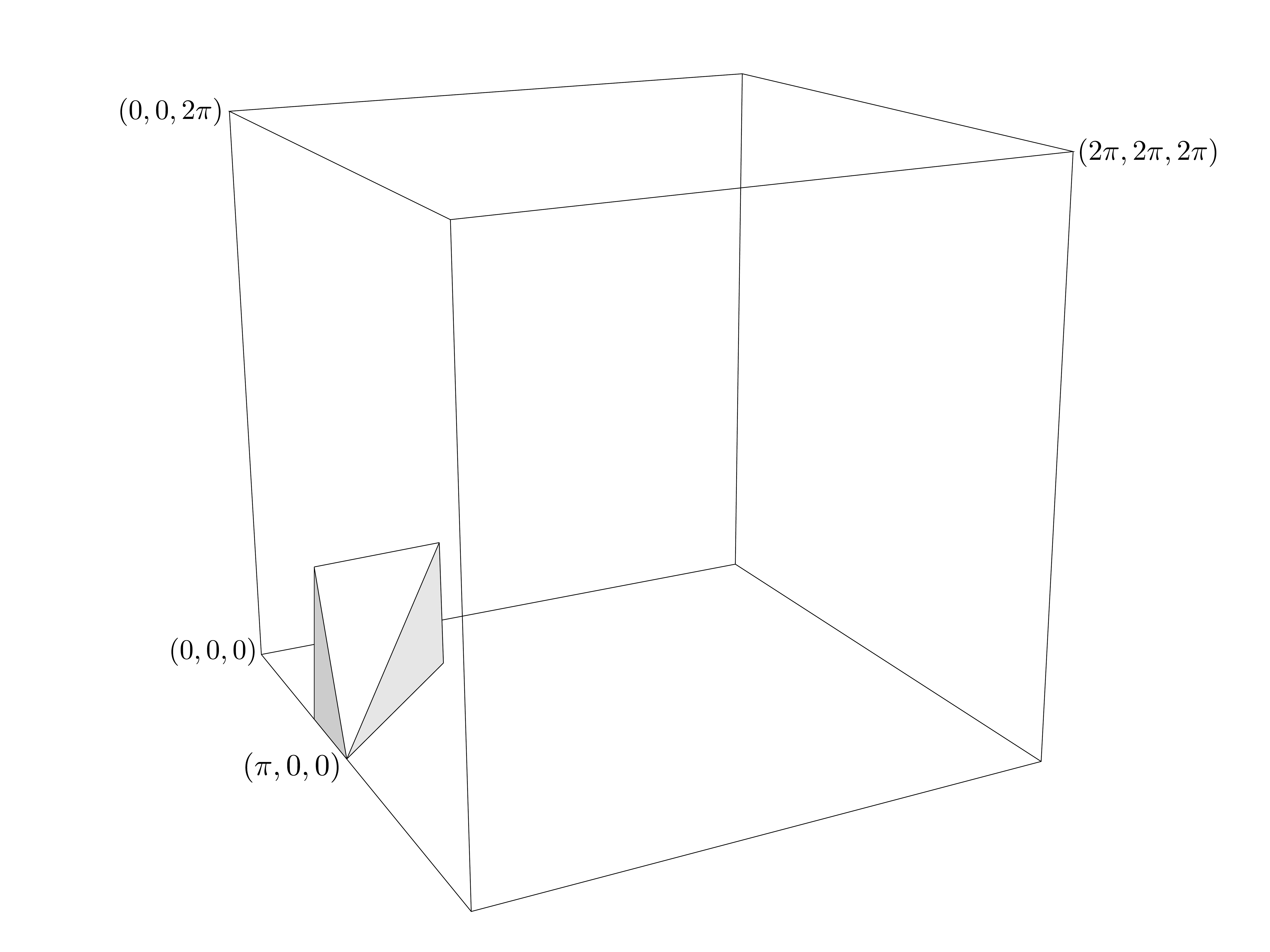}}
\caption{The figure on the left shows the portion of the moment polytope $P_6$ where $d_1>d_2,d_3$ and $(d_1)^2>(d_2)^2+(d_3)^2$. The figure on the right shows the portion of cube $[0,2\pi]^3$ where $\theta_1\in(\frac{\pi}{2},\pi)$, $\theta_2,\theta_3\in(0,\frac{\pi}{2})$, $\theta_1+\theta_2<\pi$, and $\theta_1+\theta_3<\pi$.}
\label{torussmallest}
\end{figure}

\begin{thm}
The probability that an equilateral hexagon is knotted is at most $\frac{14-3\pi}{192}$.
\end{thm}
\startproof
Let $H \in Equ_0(6)$. We will choose the $T_{135}$ triangulation of $H$ to form our set of action-angle coordinates: $\alpha:P_6\times T^{3}\mapsto Pol_0(6)$. Since almost all of $Pol_0(6)$ is a toric symplectic manifold, Theorem \ref{DH} holds for integrals over this space. First we will calculate the expected value for $H=(d_1,d_2,d_3,\theta_1,\theta_2, \theta_3)$ to have $curl=1$. Suppose that $H$ is in general position so that the lengths of the diagonals are distinct. Without loss of generality, we assume that $d_1$ is the largest of the three diagonals.  The moment polytope, $P_6$, corresponding to the $T_{135}$ triangulation is $\frac{1}{2}$ of the cube $[0,2]^3$. Therefore the volume of $P_6$ is $4$. The region where one diagonal is greater than the other two divides the moment polytope into $3$ regions with equal volume of $\frac{4}{3}$. From Lemma \ref{polytopeportionR+} and Lemma \ref{polytopeportionL+} if $d_1>d_2,d_3$ and $curl(H)=1$, then $\theta_1\in(0,\pi)$, $\theta_2\in(0,\frac{\pi}{2})$, and $\theta_3\in(0,\frac{\pi}{2})$. Additionally, if $(d_1)^2>(d_2)^2+(d_3)^2$, then $\theta_1\in(\frac{\pi}{2},\pi)$. From Lemma \ref{angles} and Lemma \ref{angles(-1,1)}, we know that $\theta_1+\theta_2<\pi$ and $\theta_1+\theta_3<\pi$, shown in Figure \ref{torussmallest}

Using standard integration, we calculate that the volume of the portion of $P_6$ where $d_1>d_2$, $d_1>d_3$, and $(d_1)^2>(d_2)^2+(d_3)^2$ is equal to $\frac{2(\pi-2)}{3}$. The ratio of this volume out of the third of $P_6$ is $\frac{\pi}{2}-1$. The region where $\theta_1\in(\frac{\pi}{2},\pi)$, $\theta_2\in(0,\frac{\pi}{2})$, $\theta_3\in(0,\frac{\pi}{2})$, $\theta_1+\theta_2<\pi$ and $\theta_1+\theta_3<\pi$ is $\frac{1}{192}$ of the cube $[0,2\pi]^3$.

\begin{figure}[h]
\centerline{\includegraphics[width=.5\textwidth]{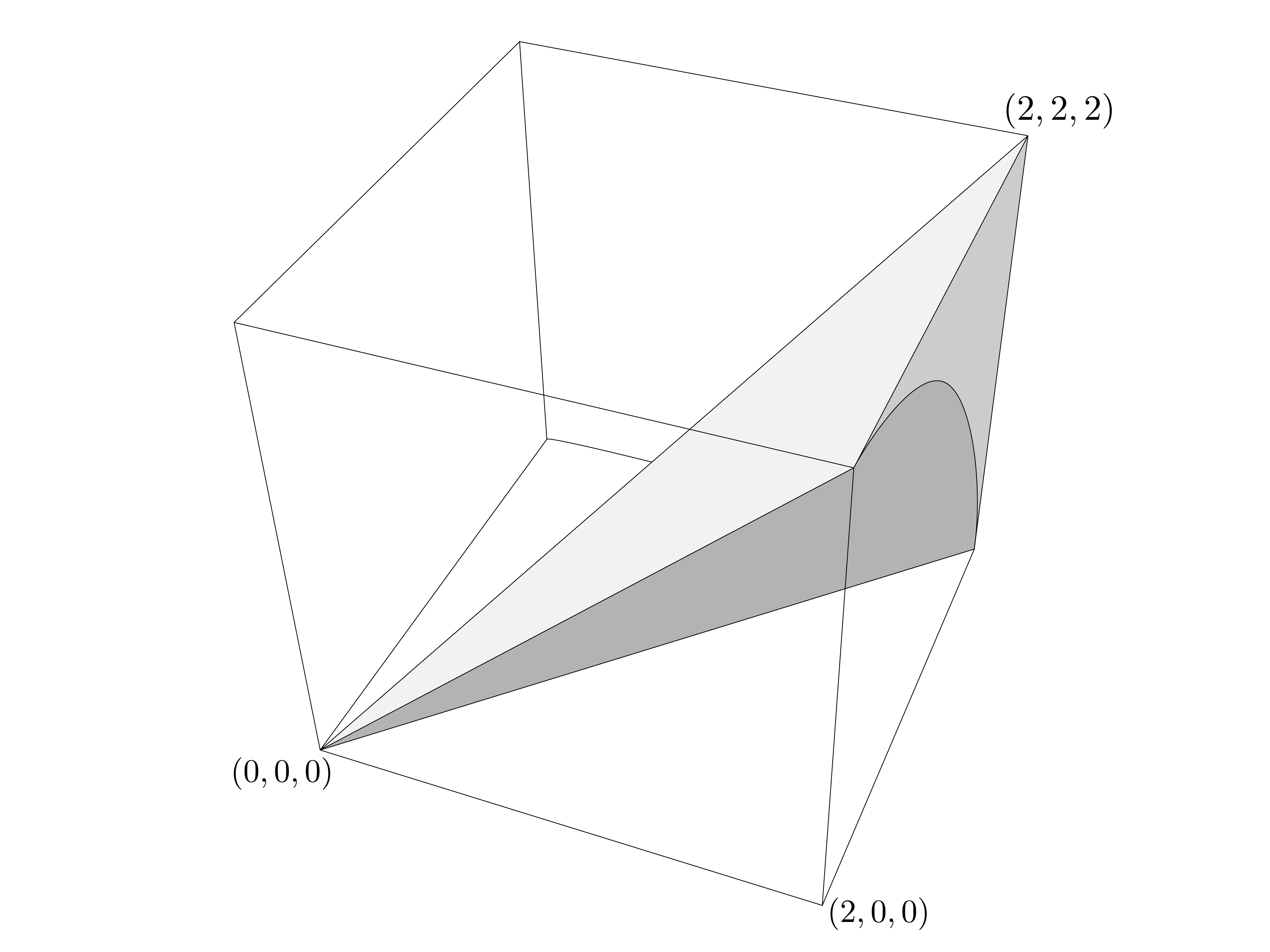}\includegraphics[width=.45\textwidth]{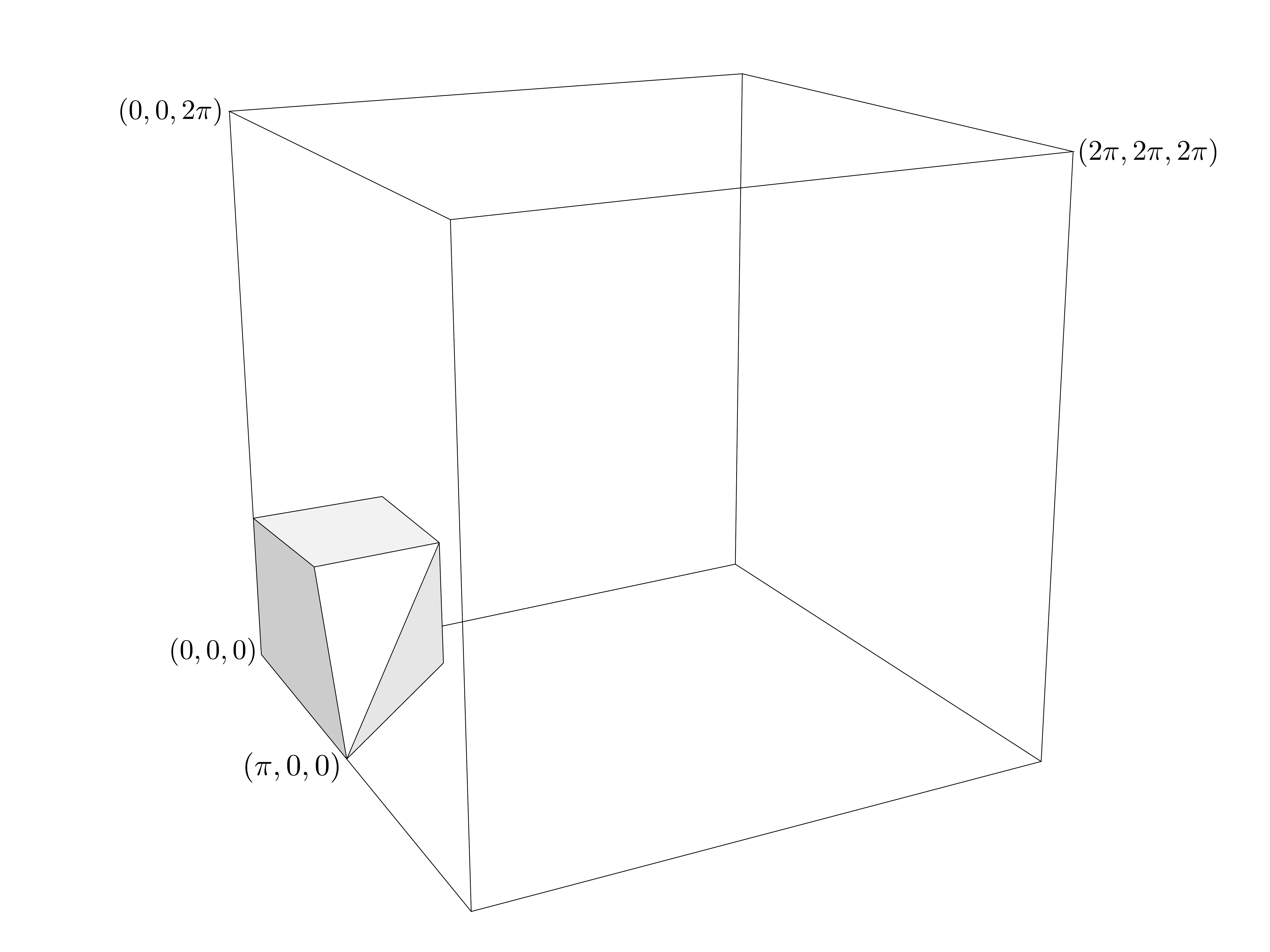}}
\caption{The figure on the left shows the portion of the moment polytope $P_6$ where $d_1>d_2,d_3$ and $(d_1)^2<(d_2)^2+(d_3)^2$. The figure on the right shows the portion of cube $[0,2\pi]^3$ where $\theta_1\in(0,\pi)$, $\theta_2,\theta_3\in(0,\frac{\pi}{2})$, $\theta_1+\theta_2<\pi$, and $\theta_1+\theta_3<\pi$.}
\label{torussmall}
\end{figure}

The portion of $P_6$ where $d_1>d_2,d_3$ and $(d_1)^2<(d_2)^2+(d_3)^2$ is $2-\frac{\pi}{2}$ of the volume of $\frac{1}{3}$ of $P_6$. The region where $\theta_1\in(0,\pi)$, $\theta_2\in(0,\frac{\pi}{2})$, $\theta_3\in(0,\frac{\pi}{2})$, $\theta_1+\theta_2<\pi$ and $\theta_1+\theta_3<\pi$, shown in Figure \ref{torussmall}, is $\frac{1}{48}$ of the cube $[0,2\pi]^3$.

 Therefore the expected value for $curl(H)=1$ is bounded above by 
$$\Big(\frac{\pi}{2}-1\Big)\Big(\frac{1}{192}\Big)+\Big(2-\frac{\pi}{2}\Big)\Big(\frac{1}{48}\Big)=\frac{7}{192}-\frac{\pi}{128}.$$
Making a similar argument for $curl(H)=-1$, we see that the knot probability is at most $$2(\frac{7}{192}-\frac{\pi}{128})=\frac{14-3\pi}{192},$$ as desired.
\finishproof \\

\newpage
\section{Acknowledgements}
The author would like to thank Ken Millett for his guidance on this project while at the University of California, Santa Barbara. The author would also like to thank Jorge Calvo, Jason Cantarella, and Clay Shonkwiler, whose work inspired this project.





\bibliographystyle{JHEP3}
\bibliography{KnottingProbability-Hake.bib}

\end{document}